\documentclass[10pt,onecolumn,aps,prd,floats,floatfix,superscriptaddress,nofootinbib]{revtex4}
\usepackage{amsmath,amsfonts,dsfont}
\usepackage{mathtools}
\usepackage{multirow}
\usepackage{graphicx}
\usepackage{color}
\usepackage{epsfig}
\usepackage{verbatim}
\begin{document}
\title{Near-Stability of a Quasi-Minimal Surface Indicated Through a \\ Tested  Curvature Algorithm}
\author{Daud Ahmad}
\thanks{daud.math@pu.edu.pk}
 \thanks{D. A. acknowledges the financial support of Punjab University, No.D/34/Est.1 Project 2013-14.}
\affiliation {Department of Mathematics, University of the Punjab, Lahore, Pakistan}
\author{Bilal Masud}
\thanks{bilalmasud.chep@pu.edu.pk}
\thanks{B. M. acknowledges the financial support of Punjab University, Sr.108 PU Project 2013-14.}
\affiliation {Center for High Energy Physics, University of the Punjab, Lahore, Pakistan}
\begin{abstract}
We decrease the $rms$ mean curvature and area of a variable surface with a fixed boundary by iterating a few times through a curvature-based variational algorithm. For a boundary with a known minimal surface, starting with a deliberately chosen non-minimal surface, we achieve up to 65 percent of the total possible decrease in area. When we apply our algorithm to a bilinear interpolant bounded by four \emph{non-coplanar} straight lines,  the area decrease by the same algorithm is only $0.116179$ percent of the original value. This relative stability suggests that the bilinear interpolant is already  a quasi-minimal surface.
\end{abstract}

\maketitle

\section{Introduction}\label{intro}
In a variational problem, we write the form (termed \emph{ansatz}) of a quantity to be minimized in such a way that what remains to be found is the value(s) of the variational parameter(s) introduced in the ansatz. An important application of the variational approach is to find characteristics of a surface (called a \emph{minimal surface}) locally minimizing its area  for a known boundary.

Minimal surfaces initially arose as surfaces of minimal surface area subject to some boundary conditions, a problem termed \textit{Plateau's problem}\cite{Osserman1986,Nitsche1989} in variational calculus. Initial non-trivial examples of minimal surfaces, namely the catenoid and helicoid, were found by Meusnier in 1776. Later in 1849, Joseph Plateau showed that minimal surfaces can be produced by dipping a wire frame with certain closed boundaries into a liquid detergent. The problem of finding minimal surfaces attracted mathematicians like Schwarz \cite{Schwarz} (who investigated triply periodic surfaces with emphasis on the surfaces called D (diamond), P (primitive), H (hexagonal), T (tetragonal) and CLP (crossed layers of parallels)), Riemann \cite{Osserman1986}, Weierstrass \cite{Osserman1986}  and R. Garnier \cite{Garnier}. Later, significant results were obtained by  L. Tonelli \cite{Tonelli},  R. Courant \cite{Courant1} \cite{Courant2}, C. B. Morrey \cite{Morrey1} \cite{Morrey2}, E. M. McShane \cite{Shane}, M. Shiffman \cite{Shiffman}, M. Morse \cite{MT}, T. Tompkins \cite{MT}, Osserman \cite{Osserman1970}, Gulliver \cite{Gulliver}, Karcher \cite{Karcher} and others.

The characteristic of a minimal surface resulting from a vanishing variation in its area is that its mean curvature should be zero \emph{throughout}; see, for example, section 3.5 of ref. \cite{docarmo} for the partial differential equation (PDE) $H=0$ this demand implies. This results in a prescription that to find  a minimal surface we should find a surface for which the mean curvature function is zero at each point of the surface. Numerically, we have to find a zero for each set of values, which means we have to solve a large number of problems on a large grid. In converting this problem into a variational form we minimize the mean square  \emph{functional} of the mean curvature numerator i.e. $\mu_{n}^{2}$ (given below by eq.~\eqref{rmsn}) with respect to our variational parameter. This functional is more convenient than minimizing directly the area functional which  has a square root in the integrand, but has the same extremal (i.e. a surface where this alternative functional is zero) as that of the area. This alternative can be compared with Douglas' suggestion of minimizing the Dirichlet integral that has the same extremal as the area functional. (A list of other possibilities of such functionals can be found in refs.\cite{MUComparative,Chen2009}.) An additional advantage is that for $\mu_{n}^{2}$ we know the value (that is, zero) to be achieved; for the area integral we do not know before calculations the target minimum value.

In an earlier work \cite{vms}, for a boundary composed of four straight lines, we reduced the $rms$  curvature and area of an initial surface chosen to be a bilinear interpolant. In this paper we took the ansatz for change in the surface to be proportional to the numerator of the mean curvature function for the surface; other factors in the change were the variational parameter, a function of the surface parameters (not to be confused with the variational parameters) whose form vanishes at the boundary and a vector assuring that we add a 3-vector to the original surface in 3-dimensions. In this paper we make an iterative use of the same ansatz, written in eq.~\eqref{VSN} below, to improve on this initial surface along with two other surfaces. Now we calculate afresh our variational parameter $t$ for each iteration $n$. An important feature of our work is that at each iteration, the function we minimize (that is $\mu_{n}^{2}$ of eq.~\eqref{rmsn}),  remains a polynomial in our variational parameter $t$. We used a computer algebra system to carry out our calculations and thus computational problems stopped us from actually implementing all the iterations. We somewhat avoided this impasse for a simpler case using a curve rather than a surface. This replaced the target flat surface by a straight line and the mean curvature by a second derivative. In this case too we had to eventually minimize a polynomial. In each case, the resulting value of the variational parameter specifies the curve or surface at each iteration.

In this paper, before reducing the area of the bilinear interpolant for which no minimal surface is completely known, we first considered two cases (a hemiellipsoid bounded by  an elliptic curve and a hump-like surface spanned by four straight lines) where we already know a minimal surface for the boundary but we deliberately start with a non-minimal surface for the same boundary.  We did this to find what fraction of the total  possible area decrease (area of the starting non-minimal surface minus that of the known minimal surface) we achieve in the computationally manageable few iterations.
If the target surface is flat, achieving it is sufficient but not necessary. This is because, at least according to the traditional definition of a minimal surface as a surface with zero mean curvature, for a fixed boundary we may have more than one minimal surface;  vanishing mean curvature is a solution of an equation obtained by setting to zero the derivative of the surface with respect to the variational parameter in its modified expression and thus a solution surface can be guaranteed to be  \emph{only a locally minimal surface} in the set of all the surfaces generated by this modification and hence is not unique. Thus another question about our algorithm is whether or not it gets stuck in some other (locally) "minimal" curve or surface before it reaches the straight line or a flat surface. We tested our ansatz for two cases with known minimal surfaces; of course for neither case the initial surface we chose was a  minimal one.

After testing in this way our algorithm, we used it to find a surface of smaller $rms$ mean curvature and area for a boundary for which no minimal surface is known. For such a case we have tried to judge if an initial surface chosen with a fixed boundary is stable or quasi-stable against the otherwise decreasing areas that our  algorithm can generate and took the resulting near stability as an indication that our initial surface is  a quasi-minimal surface for our boundary. In this paper, we report area reductions in a number of iterations to strengthen our premise. We have compared our previous ansatz with a simpler alternative (see eq.~\eqref{disposableiteration} below) to point out that the ansatz we used can be repeatedly used, which is not the case with each possible ansatz, and thus is a non-trivial feature of the ansatz we used and are further using in this paper.  The ansatz gives a significant reduction in areas  in case of hemiellipsoid and hump-like surfaces (with already known  minimal area) in computationally manageable few iterations whereas for the bilinear interpolant (a quasi-minimal surface of unknown minimal area) the decrease remains less than 0.1 percent of the original area. This also enabled us to  numerically work out  differential geometry related quantities for these surfaces. We have not been able to achieve a minimal surface within our  computationally manageable resources. But, considering the possibility of implementing the same broad algorithm differently (maybe with a better computer program or with less reliance on the computationally demanding algebra system), our present work may be a test-bed for improved detailed algorithms aimed at computing minimal surfaces.

Our work can be compared to others \cite{Monterde2004,MonterdeUgail2004} who have converted Dirichlet integrals into a system of linear equations which can be solved \cite{Monterde2004,Chen2009} to obtain extremals of a Dirichlet integral and thus surfaces of reduced area of a class of B\'{e}zier surfaces \cite{FD}. \cite{Chen2009} employs the Dirichlet method and the extended blending energy method to obtain an approximate solution of the Plateau- B\'{e}zier problem by introducing  a parameter $\lambda$ in the Extended Dirichlet Functional (compare eq.~(4) of the ref. \cite{Chen2009}  with our eq.~\eqref{rmsn}).  Determining this parameter $\lambda$  gives all the inner control points obtained directly as the solution of a system of linear equation. In our case, determining the variational parameter $t$  gives us variationally improved surfaces. The constraint that the mean curvature is zero is too strong and in most of the cases   there is no known B\'{e}zier surface \cite{Chen2009} with zero mean curvature.

Xu et al.~\cite{Xu2009} study approximate developable surfaces  and approximate minimal surfaces (defined as the minimum of the norm of mean curvature) and obtain tensor product B\'{e}zier surfaces  using a nonlinear optimization algorithm. Hao et al. \cite{HaoYong2013} find the parametric surface of minimal area defined on a rectangular parameter domain among all the surfaces with prescribed borders using an approximation   based on Multi Resolution Method using B-splines. Xu and Wang \cite{Xu2010} study quintic parametric polynomial minimal surfaces and their properties. Pan and Xu \cite{QuingPan2011} construct minimal subdivision surfaces with given boundaries using the mean curvature flow, a second order geometric PDE,  which is solved by a finite element method. For some possible applications of minimal or quasi-minimal surfaces spanning bilinear interpolants and for a literature survey, one can read the introductory section of our work \cite{vms} on \textit{Variational Minimization of String Rearrangement Surfaces} and   \cite{dabm2013} on \textit{Coons Patch Spanning a Finite Number of Curves}.

This paper is organized as follows: In  section \ref{basics}  we first point out the variational aspect of the proof of the connection between area reduction and vanishing mean curvature to provide motivation for our  ansatz. We give related definitions and the variational algorithm  based on our ansatz containing the numerator of the mean curvature function to reduce the area of a surface spanned by a fixed boundary (including the boundary composed of a finite number of given curves). In section~\ref{hessfal}, we first present a one-dimensional analogue of our ansatz of eq.~\eqref{VSN} which is applied to a curve of given length. For this simplest case, we report implementing eight iterations to get curves of reduced lengths. This is followed by applying our technique to  a hemiellipsoid surface (with an ellipse for a boundary) and  a hump-like surface (spanned by four \emph{coplanar} straight lines) to test the algorithm for the number of iterations it takes to return a surface of significantly reduced area for a boundary for which  a minimal surface is known; in these complicated cases we could implement fewer iterations. The  technique is then applied to  the bilinear interpolant bounded by four \emph{non-coplanar} straight lines, for which a minimal surface is not known. A comparison of root mean square of mean curvature and Gaussian curvature in the three  cases is also provided for further analysis of certain properties of these surfaces.  Based on this comparison,  results and remarks are presented in the final section \ref{conclusion}.

\section{Motivation for the Ansatz} \label{basics}
Our goal is to decrease the  area functional \cite{docarmo,Goetz}
\begin{equation} \label{af}
A(\mathbf x)=\int\int_{D}\left| \mathbf x_{u}(u,v)\times \mathbf x_{v}(u,v) \right|du dv,
\end{equation}
of a locally parameterized surface $\mathbf x= \mathbf x(u,v)$. Here $ D \subset R^{2} $ is a domain over which the surface $ \mathbf x(u,v) $ is defined as a map, with the boundary curve given by  $\mathbf x(\partial D)=\Gamma$  for $ 0\leq u,v\leq1 $.   $\mathbf x_{u}(u,v)$ and $\mathbf x_{v}(u,v)$ are the  partial derivatives of $ \mathbf x(u,v) $ with respect to $u$ and $v$.  The normal variation of $\mathbf{x}\left( {\bar{D}} \right)$ ($\bar{D}$ is the union of the domain $D$ with its boundary $\partial D$), determined by a differentiable function $h:\bar{D}\rightarrow R$,
is the map $\phi :\bar{D}\times \left( -\varepsilon ,\varepsilon\right)\rightarrow {{R}^{3}}$ defined by
\begin{equation}
\phi \left( u,v,t \right)=\mathbf{x}\left( u,v \right)+t\,h\left( u,v \right)\, \mathbf{N}\left( u,v \right),\text{   }\left( u,v \right)\in \bar{D},\text{   }t\in \left( -\varepsilon ,\varepsilon  \right). \label{general variation}
\end{equation}
For each fixed $t\in \left( -\varepsilon ,\varepsilon  \right)$, the map ${{\mathbf{x}}^{t}}:D \rightarrow {{R}^{3}}$ given by ${{\mathbf{x}}^{t}}(u,v)=\phi (u,v,t)$ is a parameterized surface. Let $E,F,G$ and ${{E}^{t}},{{F}^{t}},{{G}^{t}}$ denote the first fundamental magnitudes \cite{docarmo} of $\mathbf{x}(u,v)$ and ${{\mathbf{x}}^{t}}(u,v,t)$, $e,f,g$ the second fundamental magnitudes of $\mathbf{x}(u,v)$, $\mathbf{N}=\mathbf{N}(u,v)$ the unit normal to $\mathbf{x}(u,v)$ and $H$, the  mean curvature function of the surface $\mathbf{x}(u,v)$. Ref.~\cite{docarmo} shows that the derivative of the area integral
\begin{equation}
  A(t)=\int\limits_{{\bar{D}}}{\sqrt{{{E}^{t}}{{G}^{t}}-{{({{F}^{t}})}^{2}}}}\text{ }dudv,
\end{equation}
at $t=0$ is
\begin{equation}\label{areaderivative}
  {A}'(0)=-\int\limits_{{\bar{D}}}{2hH\sqrt{EG-{{F}^{2}}}}dudv,
\end{equation}
with
\begin{equation}
H(u,v) =\frac{G e -2 F f+ E g}{2(E G -F^{2})}. \label{mc}
\end{equation}

For our chosen initial non-minimal surfaces, $H(u,v)$ is non-zero. The basic idea of our algorithm is that it is always possible to decrease the area and thus to get a negative value of ${A}'(0)$ by choosing the differential function $h(u,v)$ to be proportional to the mean curvature function $H(u,v)$. Because our target is only a sign of the product $h(u,v)H(u,v)$, we simplify our work by using only the numerator of the mean curvature $H$ given by \eqref{mc}. (This is also done in ref. \cite{BCDH} following ref. \cite{Osserman1986}  that   ``for a locally parameterized surface, the mean curvature vanishes when the numerator  of the mean curvature is equal to zero".) Thus our iterative scheme for the successive surfaces $\mathbf x_{n}(u,v) \, (n=0,1,2,...)$ is
\begin{equation}\label{VSN}
\mathbf x_{n+1}(u,v,t)=\mathbf x_{n}(u,v)+t \, m_{n}(u,v) \,
\mathbf{N}_{n},
\end{equation}
where $t$ is our variational parameter and
\begin{equation}\label{vpn}
m_{n}(u,v)=b(u,v) \, H_{n},
\end{equation}
with $b(u,v)$  chosen so that the variation at the boundary curves  is zero. $H_{n}\, (\text{for}\, n=0,1,2,...)$ denotes the numerator of the mean curvature function eq.~\eqref{mc} of the non-minimal surface  $\mathbf x_n(u,v)$ and is given by
 \begin{equation}\label{mcn}
 H_{n}=e_{n} \, G _{n}-2F_{n} \,  f_{n}+g_{n} \,  E_{n}.
\end{equation}
Our ansatz in eq.~\eqref{VSN} can be compared with the arbitrary variation of eq.~\eqref{general variation} to see the choices made in our ansatz. For $\mathbf{x}_n(u,v)$, we denote  by $E_{n}, F_{n}, G_{n}, e_{n}, f_{n} $ and $g_{n}$  the fundamental magnitudes and by $\mathbf{N}_{n}(u,v)$ the numerator of the unit normal to the surface $\mathbf {x}_n (u,v)$. The subscript $n$ is used not only to denote the numerator of the quantities but also to denote the $nth$ iteration. For non-zero $n$ each of the above functions of an iterative surface is also a  function of the corresponding $t$ in addition to usual dependence on the parameters $u$ and $v$ of the surface. The functional dependence on $t$ is always a polynomial one. That is,
\begin{equation}\label{MCN}
H_{n}=H_{n}(u,v,t) = E_{n} \, g_{n}-2F_{n}  \, f_{n}+G_{n}  \, e_{n} = \sum^{6}_{i=0} (p_{i} (u,v)) \hspace{0.1cm} t^{i}.
\end{equation}
 We have written in the introduction that in place of the problematic area functional, what we minimize to find $t$ is
\begin{equation}\label{rmsn}
\mu _{n}^{2}(t)=\int_{0}^{1}{\int_{0}^{1}{H_{n}^{2}}}(u,v,t) \, dudv=\sum\limits_{i=0}^{m}{{{t}^{j}}\left( \int_{0}^{1}{\int_{0}^{1}{{{q}_{j}}\left( u,v \right)dudv}} \right)}.
\end{equation}
Because of eq.~\eqref{MCN}, $H_{n} ^{2} (u,v,t)$ in our expression for $\mu _{n}^{2}(t)$ is also a polynomial in  $t$ with real coefficients of $t^j$ for $j=0,1,2,...,10$ that we call $q_{j} (u,v)$; there are no powers of $t$ higher than 10 in the polynomials as can be seen from the expressions for $E_{n}(u,v,t)$, $F_{n}(u,v,t)$ and $G_{n}(u,v,t)$ which are quadratic in $t$ and $e_{n}(u,v,t)$, $f_{n}(u,v,t)$ and $g_{n}(u,v,t)$ which are cubic in $t$. Integrating (numerically if necessary) these coefficients with respect to $u$ and $v$ in the range $0\leq u,v\leq 1$ we get $\mu _{n}^{2}(t)$ that we minimize with respect to $t$. The resulting value of $t$ completely specifies a $new$ surface $\mathbf x_{n+1}(u,v)$. For this value $t_{min}$ of $t$, the new surface $\mathbf x_{n+1}(u,v)$ has less $ms$  mean curvature than the $ms$  mean curvature
\begin{equation} \label{rmsntmin}
\mu_{n}^{2}= \int^{1}_{0} \int^{1}_{0}H_{n}^{2}(u,v,t=t_{min})\hspace{0.2cm}dudv
\end{equation}
of the surface $\mathbf x_{n}(u,v)$. This surface is also expected to have less area and in our actual calculations we found that as we decrease our alternative functional $\mu_{n}^{2}$, the area functional of the surface spanning our fixed boundary also decreases.   The $rms$ mean curvature of both our starting surface and the one achieved after one variational area reduction remains non-zero and we re-use our variational algorithm for the resulting surface a number of times. Mean curvature tells how much the two principal curvatures~\cite{docarmo} of the surface cancel. To get an estimate of the absolute sizes of the principal curvatures we also calculated for each iteration the mean square ($ms$)  of Gaussian curvature ~\cite{docarmo} for $ 0\leq u\leq1 $ and $ 0\leq v\leq1 $. This we call $\nu_{n}$, given by the following expression:
\begin{equation} \label{rgcn}
\nu_{n}^{2} (t)=\int^{1}_{0} \int^{1}_{0} K_{n} ^{2} (u,v,t = t_{min}) \hspace{0.2cm} dudv,
\end{equation}
where $K_{n}$ is the numerator of the Gaussian curvature $K$. Gaussian curvature is the product of principal curvatures and mean curvature is average of the principal curvatures. Thus a ratio of Gaussian curvature and square of mean curvature is dimensionless. Using  eqs.~\eqref{rmsn} and \eqref{rgcn}, this ratio is
\begin{equation}\label{ratiogcmc0}
    \frac{\nu_{n}}{\mu_{n}^{2}} = \frac{\left({\int^{1}_{0} \int^{1}_{0}K_{n}^{2}(u,v,t = t_{min})\hspace{0.2cm}dudv}\right)^{1/2}}{\int^{1}_{0} \int^{1}_{0}H_{n}^{2}(u,v,t = t_{min}) \hspace{0.2cm}dudv}.
\end{equation}
Now for $\mathbf{x}_{n}(u,v)$, we compute $E_{n}(u,v,t=t_{min})$, $F_{n}(u,v,t=t_{min})$ and $G_{n}(u,v,t=t_{min})$ and denote the area integral eq.~ \eqref{af} as $A_{n}$ which is given by the following expression
\begin{equation} \label{arean}
A_{n}= \int^{1}_{0} \int^{1}_{0} \sqrt{E_{n} G_{n}-F_{n}^{2}} \hspace{0.2cm} dudv.
\end{equation}
The ansatz eq.~\eqref{VSN}, for the first order reduction in the area  of a non-minimal surface $\mathbf x_{0}(u,v)$ is
\begin{equation}\label{VS1}
\mathbf x_{1}(u,v,t)=\mathbf x_{0}(u,v)+t\,m_{0}(u,v) \mathbf{N}_{0},
\end{equation}
where $t$ is our variational parameter. Here $m_{0}(u,v), H_{0}(u,v)$ are given by eqs.~ \eqref{vpn}, \eqref{MCN} and  $\mathbf{N}_0 (u,v)$  is the   unit normal  to the non-minimal surface $\mathbf{x}_0(u,v)$ for $n=0$. For $n=0$, the ratio \eqref{ratiogcmc0} of the $rms$ of Gaussian curvature to the $ms$ of mean curvature is given by
\begin{equation}\label{ratiogcmc0}
  \frac{\nu_{0}}{\mu_{0}^2} = \frac{\left({\int^{1}_{0} \int^{1}_{0}K_{0}^{2}\hspace{0.2cm}dudv}\right)^{1/2}}{\int^{1}_{0} \int^{1}_{0}H_{0}^{2}\hspace{0.2cm}dudv},
\end{equation}
where $K_{0}$ is the numerator of the Gaussian curvature $K$. With the above notation,  eq.~ \eqref{af} becomes, for $\mathbf x_0 (u,v)$,
\begin{equation} \label{area0}
A_{0}= \int^{1}_{0} \int^{1}_{0} \sqrt{E_{0} G_{0}-F_{0}^{2}} \hspace{0.2cm} dudv.
\end{equation}
For $n=1, 2, ...$ in eq.~\eqref{VSN} gives us the  surfaces $\mathbf x_{2}(u,v)$, $\mathbf x_{3}(u,v), ...$ of reduced area and  related quantities may be computed from eqs.~ \eqref{vpn} to \eqref{arean}.  In order to see a geometrically meaningful (relative) change in area we calculate the dimensionless area ratios. For the cases with known minimal surface, let $A_{0}$ be the initial area and $A_{f}$ the area of the known minimal surface
(which we took, for these ratio calculations, as the corresponding flat surface; the subscript $f$ reminds us of this). We can define the maximum possible change to be achieved as $\triangle A_{max}\, =
A_{0}-A_{f}$. Let $A_{i}$ for $i=0, 1, 2, ...$ be the area of the surface $\mathbf{x}_{i}(u,v)$ obtained through the $ith$ iteration. Then the difference of the areas in the $ith$ and $jth$ iteration, with $i<j$, is denoted by $\triangle A_{ij}$ . When we know the minimal
surface, the percentage decrease $p_{ij}$ in area can be computed by multiplying the quotient of $\triangle A_{ij}$ and $\triangle A_{max}$
by $100$. Thus we have
\begin{equation}\label{percentdecknown}
  p_{ij} = 100 \frac{\triangle A_{ij}}{\triangle A_{max}} = 100 \frac{A_{i}-A_{j}}{A_{0}-A_{f}} \, \hspace{0.25cm} \text{for} \,i<j \hspace{0.25cm} \text{and} \, i,j=0, 1, 2, ....
\end{equation}
When we do not know the minimal surface, the percentage decrease $q_{ij}$ in area can be computed by using
\begin{equation}\label{percentdecunknown}
  q_{ij} = 100 \frac{A_{i}-A_{j}}{A_{0}} \, \hspace{0.25cm} \text{for} \, i<j \hspace{0.25cm} \text{and} \, i,j=0, 1, 2, ....
\end{equation}
\section{The Efficiency Analysis of the Ansatz eq.~ \eqref{VSN}} \label{hessfal}
\begin{figure}[htb!]
\begin{center}
\includegraphics[width=60mm]{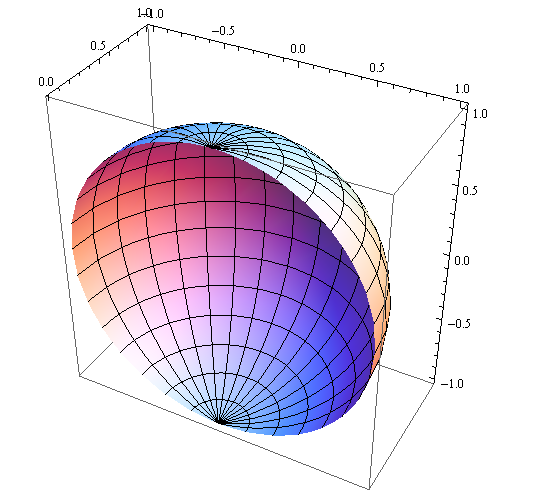}
\end{center}
\caption{Hemiellipsoid spanned by an elliptic curve for $0\le u\le \pi $  and $0\le v\le \pi$.} \label{hemiellipsoidgraphinitial}
\end{figure}

\begin{figure}[htb!]
\begin{center}
\includegraphics[width=60mm]{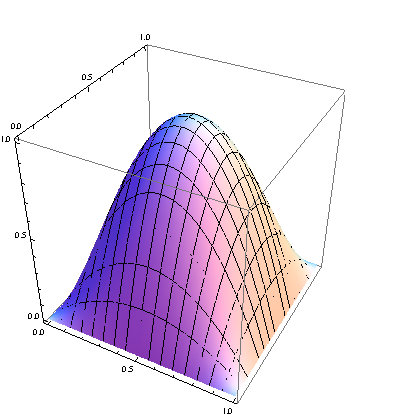}
\end{center}
\caption{A hump-like surface spanned by four coplanar straight lines with the parametrization  $\mathbf x_{0} \left(u,v \right) = (u, v, 16u v (1-u) (1-v))$ for $0\leq u,v \leq 1$.} \label{surface1}
\end{figure}

\begin{figure}[htb!]
\begin{center}
\includegraphics[width=60mm]{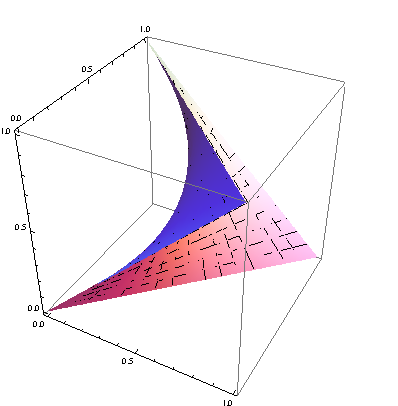}
\end{center}
\caption{A bilinear interpolant spanned by four non-coplanar straight lines with the parametrization $\mathbf x_{0} \left(u,v\right) = (r(u+v-2uv), v, u)$ for $r=1$ and  $0\leq u,v\leq 1$.} \label{binterpolation1}
\end{figure}

In this section we apply the technique introduced in  section  \ref{basics}  to reduce the area of  a variety of non-minimal surfaces. We start by reporting calculations for a hemiellipsoid bounded by  an elliptic curve for which a minimal surface is the elliptic disc. This starting non-minimal surface is given by  the parametrization
\begin{equation}\label{hemiellipsoid}
  \mathbf{x}_{0}(u,v)=(\sin u\cos v,b\sin u\sin v,c\cos u),
\end{equation}
where  $b$,  $c$ are  constants and $0\le u\le \pi $  and $0\le v\le \pi$ as shown in Fig.~\ref{hemiellipsoidgraphinitial}. For a second case, we took the starting surface as a hump-like surface spanned by four arbitrary straight lines with parametrization
  \begin{equation}\label{eqsurface1}
    \mathbf x_{0} \left(u,v \right) = (u, v, 16u v (1-u) (1-v))
\end{equation}
where $0 \leq u,v \leq 1$ as shown in Fig.~ \ref{surface1}. The target minimal surface in this case is also a flat surface, but the boundary is  instead a square.  The minimum area in this second case is that of a square of unit length.  Our success  (reported below) in significantly reducing the area in these known cases suggests it is interesting to find and analyze  area reductions in surfaces for which no minimal surface is known. Linear interpolant in the one-dimensional case (a straight line) is already minimal. But the mean curvature function for the bilinear interpolant spanning a non-planar boundary is not zero, meaning the derivative of the area function with respect to a variational parameter is non-zero.   The curvature-based algorithm we suggest manages to utilize this non-zero mean curvature and area derivative to decrease both of these quantities.  For definiteness, we took the corners
  \begin{equation}
	\mathbf x(0,0)=\mathbf r_{1}, \hspace{.5cm} \mathbf x(1,1)=\mathbf r_{2}, \hspace{.5cm}\mathbf x(1,0)=\mathbf r_{\bar{3}}, \hspace{.5cm} \mathbf x(0,1)=\mathbf r_{\bar{4}},
\end{equation}
of our bilinear interpolant as
\begin{equation} \label{ruled1}
	\mathbf r_{1}=(0,0,0),  \hspace{.5cm} \mathbf r_{2}=(r,r,0),  \hspace{.5cm} \mathbf  r_{\bar{3}}=(0,r,r),  \hspace{.5cm} \mathbf r_{\bar{4}}=(r,0,r),
\end{equation}
for a real scalar $r$. For this case, the bilinear interpolant is the following bilinear mapping from $(u,v)$ to $(x,y,z)$:
\begin{equation}\label{eqbinterpolation1}
    \mathbf x_{0} \left(u,v \right) = (r(u+v-2uv), v, u).
\end{equation}
Since $r$ is the only scale in our problem, geometrical properties do not depend on the actual value of $r$ we choose; the argument is straightforward for the dimensionless ratios we report. Thus we take the simplest choice of taking $r=1$; for this choice the bilinear interpolant is shown in Fig.~\ref{binterpolation1}.
In case of a hemiellipsoid, to check if some alterations in our algorithm can be introduced, we replaced, in our ansatz for a change in surface, the mean curvature function by a constant chosen as one. An argument shows that this simple algorithm can not be iterated: with starting value of $n=0$ and a convenient choice of $\mathbf{N}_0(u,v)=\mathbf{k}=\text{constant vector}\,\text{(say)}$ and $H_n=H_0$, $\forall\, n$, the ansatz eq.~\eqref{VSN} reduces to the following expression
\begin{equation}\label{disposableiteration}
\mathbf x_{1}(u,v,t)=\mathbf x_{0}(u,v)+\widetilde{t}\,m_{0}(u,v)\mathbf{k},
\end{equation}
where $\widetilde{t}$ is now the variational parameter. Any further iteration of this algorithm would only change the value of $\widetilde{t}$ but not the form of the resulting surface.  But the value is uniquely given by our minimization. Thus for iterating our algorithm we restored the non-trivial original mean curvature numerator $H_n,\,n=1,2,..$ defined as above by  eq.~\eqref{mcn}. The resulting decreases in $rms$ mean curvature and areas from the above mentioned starting surfaces are reported in the four subsections below. A comparison can be seen in table~\ref{table2}.
\subsection{The 1-Dimensional Analogue of the Curvature Algorithm Applied to a Curve of a Given Length}
\begin{figure}[htb!]
\begin{center}
\includegraphics[width=75mm]{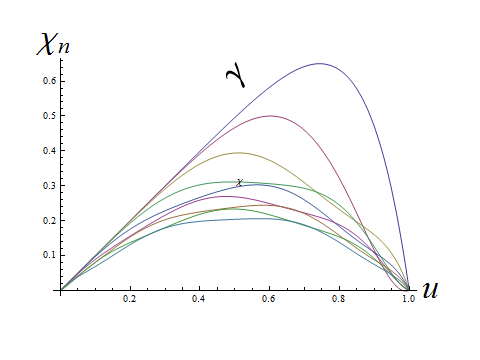}
\end{center}
\caption{Variational improvement to arc length of a curve of given length, comparing the results of eight iterations, where $\gamma$ is the initial curve and the remaining curves are the curves of reduced lengths.} \label{vimpcurveeightiterations}
\end{figure}

Before finding what our algorithm(s) yield for surfaces, we wrote a program for the variational problem of reducing arc length, targeting this time the straight line joining the end points of the curve.
For this one dimensional variational problem our ansatz for successive curves $\mathbf{\chi}_n(u), \, n=0,1,2,...$ joining the same two ends may be written in the form:
\begin{equation}\label{vc1}
  \mathbf{\chi}_{n+1}(u) = \mathbf{\chi}_n+t \, m_{n}(u)\mathcal{N}_n, \hspace{0.25cm} n=0,1,2...,
\end{equation}
where
\begin{equation}\label{vpc}
  m_n(u)=u(1-u)H_n
\end{equation}
is chosen so that it is zero at $u=0$ and $u=1$.  Here, $n$ is the iteration number. $\mathcal{N}_n$ is a unit vector in the direction normal to the curve, which practically means we take the transverse displacement from the straight line joining the two ends. In place of the
numerator of the mean curvature for the surface case, we take here $H_n$ to be  ordinary curvature which is the second derivative with respect to the curve parameter $u$. When the technique was applied to a  starting curve
\begin{equation}\label{curve0}
  \mathbf{\chi}_0(u)=(u,u-u^8),
\end{equation}
it gave  the following expressions for variationally improved curves
\begin{flalign}
  & {{\mathbf{\chi }}_{1}}\left( u \right)=\left( u,u-7.4286{{u}^{7}}+6.4286{{u}^{8}} \right), \\
 & {{\mathbf{\chi }}_{2}}\left( u \right)=\left( u,u-22.131{{u}^{6}}+40.2381{{u}^{7}}-19.1071{{u}^{8}} \right), \\
 & {{\mathbf{\chi }}_{3}}(u)=\left( u,u-40.6973{{u}^{5}}+122.159{{u}^{6}}-128.943{{u}^{7}}+46.4814{{u}^{8}} \right), \\
 & {{\mathbf{\chi }}_{4}}(u)=\left( u,u-39.6743{{u}^{4}}+177.61{{u}^{5}}-320.449{{u}^{6}}+261.908{{u}^{7}}-80.3952{{u}^{8}} \right), \\
 & {{\mathbf{\chi}}_{5}}\left(u\right)=\left(u,u-21.394{{u}^{3}}+141.344{{u}^{4}}-414.012{{u}^{5}}+605.86{{u}^{6}}-434.714{{u}^{7}}+121.916{{u}^{8}} \right), \\
 & {{\mathbf{\chi }}_{6}}\left( u \right)=\left( u,u-5.0269{{u}^{2}}+50.0547{{u}^{3}}-249.339{{u}^{4}}+622.028{{u}^{5}}-820.916{{u}^{6}}+547.646{{u}^{7}}-145.446{{u}^{8}} \right), \\
 & {{\mathbf{\chi }}_{7}}\left( u \right)=\left(u,0.6239u+6.5839{{u}^{2}}-73.1064{{u}^{3}}+327.961{{u}^{4}}-764.604{{u}^{5}}+960.763{{u}^{6}}-617.46{{u}^{7}}+159.24{{u}^{8}} \right), \\
 & {{\mathbf{\chi }}_{8}}\left( u \right)=\left( u,1.0778u-8.9886{{u}^{2}}+77.659{{u}^{3}}-334.76{{u}^{4}}+755.916{{u}^{5}}-926.532{{u}^{6}}+583.748{{u}^{7}}-148.119{{u}^{8}} \right).
\end{flalign}
Corresponding lengths of these curves are $\ell_0 = 1.7329$, $\ell_1 = 1.46525$, $\ell_2 = 1.30988$,  $\ell_3 = 1.24103$,  $\ell_4 = 1.20133$,  $\ell_5 = 1.16958$,  $\ell_6 = 1.1459$,  $\ell_7 = 1.12682$ and  $\ell_8 = 1.11081$. Percentage decreases in the lengths are denoted by \begin{equation}\label{percentdecreasecurve}
\ell _{ i  j}=100 \frac{\ell  _ i - \ell  _ j}{\ell  _0 - 1},\hspace{0.25cm} \text{where} \hspace{0.25cm}  i <  j \hspace{0.25cm} \text{and} \hspace{0.25cm}  i,  j=0,1,2....
\end{equation}
and are reported in table~\ref{curvetable}. Fig.~\ref{vimpcurveeightiterations} shows the graphs of all the curves which we could achieve before exhausting the limit of available computer resources. It can be seen from the second   column of  table~\ref{curvetable} that the length of the sequence of curves is getting closer and closer  to the shortest (unit) length joining the two ends. We also report in the next column the corresponding values of the variational parameter $t$ that gives these lengths.
\begin{center}
\begin{table}[htb!]
\caption{\textbf{Reduction in Length of a Curve of Given Length}}
\smallskip
\begin{quote}
  The decreasing lengths $\ell_i$ of the variationally improved curves  $\mathbf{\chi}_i$ along with percentage decreases in length given by $\ell_{i\,j}$($i<j$) for the corresponding optimal value $t_{min}$ of our variational parameter $t$.\end{quote}
  \vspace{3mm}
    \begin{tabular}{ | p{0.5cm} | p{1.2cm} |p{1.2cm} |p{1.2cm} | }
    \hline
   $\mathbf{\chi}_i$ & $\ell_i$ & $\ell_{i\,j} $& $t_{min}$\\ \hline
  \hspace{1.5cm}   $ \mathbf{\chi}_0 $  &  \hspace{1.5cm}  1.7329     &  \hspace{1.5cm} &  \hspace{1.5cm} \\ \hline
  \hspace{1.5cm}   $ \mathbf{\chi}_1 $  &  \hspace{1.7cm} 1.46525   & \hspace{1.5cm}  36.5206      & \hspace{1.5cm} 0.132653 \\ \hline
  \hspace{1.5cm}   $ \mathbf{\chi}_2 $  &  \hspace{1.7cm} 1.30988   & \hspace{1.5cm}   21.1996     & \hspace{1.5cm} 0.070933    \\ \hline
   \hspace{1.5cm}  $ \mathbf{\chi}_3 $  &  \hspace{1.7cm} 1.24103   & \hspace{1.5cm}   9.3933        & \hspace{1.5cm}0.061298  \\ \hline
   \hspace{1.5cm}  $ \mathbf{\chi}_4 $  &  \hspace{1.7cm} 1.20133   & \hspace{1.5cm}  5.41714      & \hspace{1.5cm}0.048743  \\  \hline
   \hspace{1.5cm}  $ \mathbf{\chi}_5 $  &  \hspace{1.7cm} 1.16958   & \hspace{1.5cm}   4.33218     & \hspace{1.5cm}0.044937  \\ \hline
   \hspace{1.5cm}  $ \mathbf{\chi}_6 $  &  \hspace{1.7cm} 1.1459     & \hspace{1.5cm}  3.2307         & \hspace{1.5cm}0.039161  \\ \hline
   \hspace{1.5cm}  $ \mathbf{\chi}_7 $  &  \hspace{1.7cm} 1.12682   & \hspace{1.5cm}  2.60332      & \hspace{1.5cm}0.037408  \\ \hline
   \hspace{1.5cm}  $ \mathbf{\chi}_8 $  &  \hspace{1.7cm} 1.11081   & \hspace{1.5cm}   2.18471     & \hspace{1.5cm}0.034467\\ \hline
   \end{tabular} \label{curvetable}
   \end{table}
\end{center}

\subsection{The Hemiellipsoid Results}
For constants $b$  and $c$ and $0\le u, v \le \pi $, fundamental magnitudes, mean curvature and the area  of the initial surface hemiellipsoid  of eq.~ \eqref{hemiellipsoid} are
\begin{align}
{{E}_{0}}=&
\begin{aligned}[t]
&{{\cos }^{2}}u\left( {{b}^{2}}{{\sin }^{2}}v+{{\cos }^{2}}v \right)+{{c}^{2}}{{\sin }^{2}}u,
\end{aligned}\\
 {{F}_{0}}=&
\begin{aligned}[t]
&\left( {{b}^{2}}-1 \right)\sin u\cos u\sin v\cos v,
\end{aligned}\\
{{G}_{0}}=&
\begin{aligned}[t]
&{{\sin }^{2}}u\left( {{b}^{2}}{{\cos }^{2}}v+{{\sin }^{2}}v \right),
\end{aligned} \\
 {{e}_{0}}=&
\begin{aligned}[t]
&-bc\sin u,
\end{aligned}\\
{{f}_{0}}=&
\begin{aligned}[t]
&0,
\end{aligned}\\
{{g}_{0}}=&
\begin{aligned}[t]
&-bc{{\sin }^{3}}u,
\end{aligned}\\
  H_{0}(u,v) =&
\begin{aligned}[t]
 &\left(-2 b\, c\,{{\sin }^{3}}u\left( \left( {{b}^{2}}-2{{c}^{2}}+1 \right)\cos (2u)+3{{b}^{2}}+2{{c}^{2}}+3 \right)-4b\left( {{b}^{2}}-1 \right)c{{\sin }^{5}}u\cos (2v) \right)/8,
 \end{aligned}\\
 A_{0} =&
\begin{aligned}[t]
 &\int_{0}^{1}{\int_{0}^{1}{\sqrt{\sin^{2} u({{b}^{2}}\cos^{2} u+{{c}^{2}}\sin^{2} u({{b}^{2}}\cos^{2} v+\sin^{2} v))}\text{ }du\text{ }dv}}.
 \end{aligned}
\end{align}
In particular for $b=1, c=1$, these quantities reduce to
$E_{0}=1,F_{0}=0,G_{0}=\sin^{2}u,e_{0}=-\sin u,f_{0}=0,g_{0}=-\sin ^{3}u$, $H_{0}=-2\sin ^{3}u, A_{0}= 6.28319$. The corresponding function $b(u,v)$ is defined as
\begin{equation}\label{buvhemiformc1}
  b(u,v) = v (\pi - v).
\end{equation}
Since the success of our algorithm depends only on the sign of the derivative of area in eq.~\eqref{areaderivative}, we can modify  the previously used changes $h(u,v)$ in our surface that do not alter the signs of the area derivative. Utilizing this freedom, we replaced the variable normal to the surface by a fixed unit vector that gives the same sign of the area derivative. For our hemiellipsoid case, $\mathbf{k}=(0,1,0)$ has this property and we replaced the normal to the surface by this simpler 3-vector. Below we give the reduction in area of the hemiellipsoid eq.~\eqref{hemiellipsoid} for the two cases, one for which  $H_0$ is replaced by $1$ and one for which $H_0$ is numerator of the mean curvature. Replacing $H_0$  by $1$ in eq.~\eqref{vpn} along with   ~\eqref{buvhemiformc1} gives the first variational surface as
\begin{equation}
  \mathbf{x}_{1}(u,v,t) = \left(\sin (u) \cos (v),t (\pi -v) v+\sin (u) \sin (v),\cos (u)\right).
\end{equation}
Thus the coefficients of $t^i$ in the expansion of  the usual numerator of the mean curvature $ H_{1}(u,v,t)$ and the mean square of  the mean curvature $\mu_{1}^{2}(t)$ (eq.~\eqref{rmsn} for $n=1$)   are given by
\begin{align}
     p_0 =&
     \begin{aligned}[t]
     &-2 \sin ^3 u,
     \end{aligned}\\
     p_{1}=&
     \begin{aligned}[t]
     & -2 \sin ^2 u \, (\sin v+2 (\pi -2 v) \cos v),
     \end{aligned}\\
      p_{2}=&
      \begin{aligned}[t]
      & -\frac{1}{4} (\pi -2 v)^2 \sin u \,  (\cos (2 (u-v))+\cos (2 (u+v))-2 \cos (2 u)+6 \cos (2 v)+6),
      \end{aligned}\\
       p_3 =&
       \begin{aligned}[t]
       & -(\pi -2 v)^3 \cos v,
       \end{aligned}\\
        \mu_{1}^{2} (t) = &
        \begin{aligned}[t]
        &1270.43 t^6+1724.78 t^5+1465.01 t^4+813.722 t^3+317.473 t^2+85.3333 t+12.337.\label{msmchemin1}
        \end{aligned}
\end{align}
Minimizing $\mu_{1}^{2} (t)$ with respect to $t$ gives us $t_{min1}=-0.351571$, so that
\begin{equation}
  \mathbf{x}_1(u,v)=(\sin (u) \cos (v),\sin (u) \sin (v)-0.351571 (\pi -v) v,\cos (u)).
\end{equation}
We found $A_{1}=4.70625$ and thus the percentage decrease in area is given by $p_{01}=50.1954$. Similarly we calculated
\begin{equation}
  \mu_{2}^2=21975.9 t^6-9141.25 t^5+5060.49 t^4-1145.39 t^3+363.123 t^2-40.3821 t+1.73308.
\end{equation}
In this case $t_{min2}=0.0706353$ and hence
\begin{equation}
  \begin{split}
    \mathbf{x}_2(u,v)  =& \, (\sin u  \cos v,\,0.0706\,(\pi -v) \,v\, (0.9888\,{{(1.5708-v)}^{2}} \sin u{{\,\cos }^{2}}u{{\,\sin }^{2}}v\,+ \\
 & \, ({{\sin }^{2}}u{{\,\sin }^{2}}v+{{(\sin u\, \cos v+0.7031\,v-1.1045)}^{2}})((1.1045-0.7031\,v) \\
 & \, \cos v-\sin u)-{{\sin }^{2}}u\text{ }(\sin u-0.7031\, \sin v+(0.7031\,v-1.1045) \\
 & \,  \cos v))+ \sin u\, \sin v-0.3515\,(\pi -v)\,v,\cos u).
  \end{split}
\end{equation}
Here, $A_{2}=4.44025$ and $p_{12}=8.4671$. The percentage decrease for the full area reduction at this stage is $p_{02}=58.6625$. Continuing the process, we find the mean square mean curvature for $n=3$.  Here
\begin{equation}
    \mu_{3}^2 = 6317.67 t^6 - 656.284 t^5+1498.48 t^4-104.093 t^3+239.168 t^2 -10.7401 t+0.401386.
\end{equation}
In this case $t_{min3} = 0.0226436$ and this  gave us $\mathbf{x}_3(u,v)$ (a lengthy expression not reproduced here) for which the area comes out to be $A_3=4.4025$, $p_{23}=1.20025$. We found $p_{03}=59.8627$. These results are presented in the table~\ref{table2}. Below we give results for the reduction in area of the hemiellipsoid eq.~\eqref{hemiellipsoid} for which we restored the non-trivial actual mean curvature numerator $H_0$. In this case we have been able to produce only two iterations.  For  $b=1, \, c=1$, $b(u,v)=v(\pi-v)$   and the expressions for $m_0(u,v)$, the first variational surface $\mathbf{x}_1(u,v,t)$ and $\mu_{1}^2 $, mean square of the mean curvature, are
\begin{align}
 {{m}_{0}}(u,v) =&
 \begin{aligned}
 &-2(\pi -v)v{{\sin }^{3}}u,
 \end{aligned} \\
  {{\mathbf{x}}_{1}}(u,v,t)=&
  \begin{aligned}
  &(\sin u\cos v,\sin u\sin v-2t\text{ }(\pi -v)\text{ }v\text{ }{{\sin }^{3}}u,\cos u),
 \end{aligned}\\
  \mu _{1}^{2}=&
  \begin{aligned}
  &39774.1\text{ }{{t}^{6}}-41607.9\text{ }{{t}^{5}}+22816.4\text{ }{{t}^{4}}-7479.4\ {{t}^{3}}+1683.3\text{ }{{t}^{2}}-219.9\text{ }t+12.3.
  \end{aligned}
\end{align}
In this case $t_{min1} = 0.148252$. Accordingly, the first variational surface is given by
\begin{equation}
\begin{split}
        \mathbf{x}_1(u,v) =  &\, ( \sin  u \cos  v, \sin u((0.6986-0.2224v) v\, {{\cos }^{2}}u+(0.0741v-0.2329)v \, {{\sin }^{2}}u+ \\& \, 0.2224 \, {{v}^{2}}-0.6986 \,v+\sin v), \cos  u).
\end{split}
\end{equation}
The percentage decrease in the $rms$ mean curvature of the numerator of the mean curvature is 74.1475. Area $A_1=4.32641$ and the percentage decrease in area $p_{01}=62.2861$, as summarised below in table~\ref{table2}. Continuing the process as above we find the dependence on our variational parameter $t$ of the mean square of  the mean curvature  is
\begin{equation}
  \mu_{2}^{2}=77354.5{{t}^{6}}+4011.85{{t}^{5}}+9417.35{{t}^{4}}-174.627{{t}^{3}}+701.03{{t}^{2}}-32.3979t+0.8246.
\end{equation}
Here $t_{min2}= 0.023$ and related results are summarized in table~\ref{table2}. It can be seen from  table~\ref{table2} that this choice of $H_0(u,v)$ gives a better reduction in the area even for the first iteration.
\subsection{A Hump-Like Surface \eqref{eqsurface1} Spanned by Four Boundary Coplanar Straight Lines}\label{hessfalhs}
\begin{figure}[htb!]
\begin{center}
\includegraphics[width=60mm]{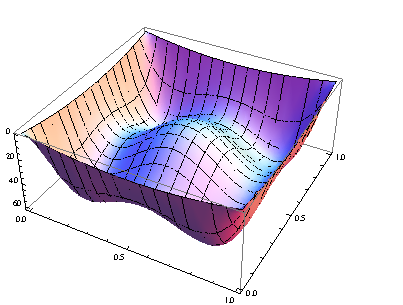}
\end{center}
\caption{$H_0(u,v)$, the numerator of the initial mean curvature of the hump-like surface $\mathbf{x}_0(u,v)$ for $0\leq u,v\leq 1$.} \label{mc0surface1}
\end{figure}

\begin{figure}[htb!]
\begin{center}
\includegraphics[width=60mm]{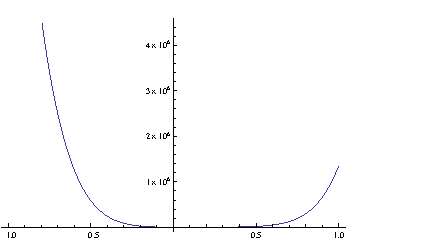}
\end{center}
\caption{Mean square $\mu_{1}^{2} (t)$ of mean curvature  of the surface $\mathbf{x_1}(u,v,t)$.} \label{graphh1nsqt}
\end{figure}

\begin{figure}[htb!]
\begin{center}
\includegraphics[width=60mm]{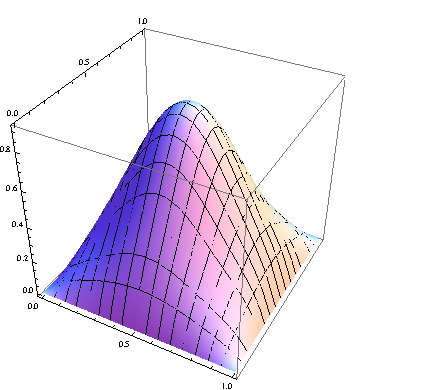}
\end{center}
\caption{The surface $\mathbf{x_1}(u,v,t)$ for $t=t_{min1} = 0.088933$ for $0\leq u,v\leq 1$.} \label{graphvarx1Mint}
\end{figure}

\begin{figure}[htb!]
\begin{center}
\includegraphics[width=60mm]{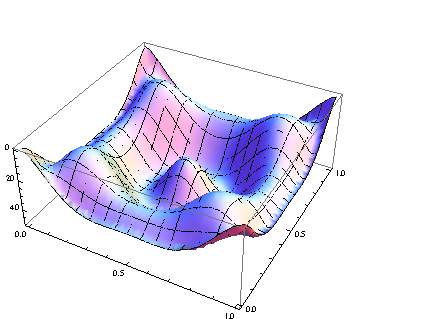}
\end{center}
\caption{Mean curvature of the surface $\mathbf{x_1}(u,v,t)$ for $t=t_{min1} = 0.088933$ for $0\leq u,v\leq 1$.} \label{graphcH1Mint}
\end{figure}

\begin{figure}[htb!]
\begin{center}
\includegraphics[width=60mm]{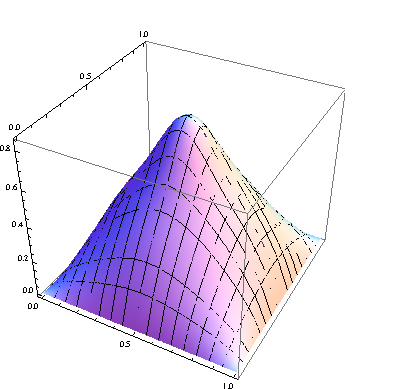}
\end{center}
\caption{Surface $\mathbf{x_2}(u,v,t)$ for $t=t_{min2} = 0.0441$ for $0\leq u,v \leq 1$.} \label{graphvarx2Mint2}
\end{figure}

\begin{figure}[htb!]
\begin{center}
\includegraphics[width=60mm]{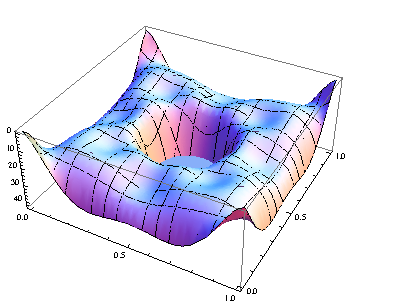}
\end{center}
\caption{Mean curvature of the surface $\mathbf{x_2}(u,v)$ for $t_{min2} = 0.0441$ for $0\leq u,v \leq 1$.} \label{graphcH2Mint2}
\end{figure}

\begin{figure}[htb!]
\begin{center}
\includegraphics[width=60mm]{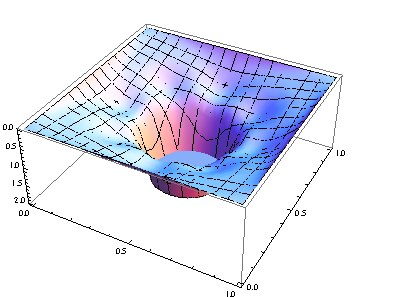}
\end{center}
\caption{$m_2(u,v)$ (eq.~\eqref{vpn} for $n=2$) for the surface $\mathbf{x_3}(u,v)$ for $0\leq u,v \leq 1$.} \label{m2graph}
\end{figure}

We apply the ansatz eq.~\eqref{VSN} along with eq.~\eqref{vpn} mentioned in section \ref{basics} to the surface $\mathbf x(u,v)$ given by eq.~ \eqref{eqsurface1} bounded by four coplanar straight lines $0\leq u,v\leq 1$. A function $b(u,v)$ whose variation at the boundary curves is zero is given by
\begin{equation}\label{humpbuv}
  b(u,v)=u v (1-u)(1-v).
\end{equation}
A convenient possible choice for the alternative to the  unit normal $\mathbf N(u,v)$ to this surface
is  $\mathbf k=(0,0,1)$ which makes a small angle with this  unit normal and thus does not change the sign of the area derivative. The surface given by eq.~\eqref{eqsurface1} is a non-minimal surface and the $xy$-plane  bounded by $0\leq u,v\leq1$ is a minimal surface spanning its boundary. The fundamental magnitudes of the initial surface eq.~\eqref{eqsurface1} are
\begin{align}
  E_0  = &
  \begin{aligned}[t]
  &1+256{{v}^{2}}{{\left( 1-2u \right)}^{2}}{{\left( 1-v \right)}^{2}},\label{humpE0}
  \end{aligned}\\
  F_0 = &
  \begin{aligned}[t]
  &256\text{ }u\text{ }v\left( 1-3u+2{{u}^{2}} \right)\left( 1-3v+2{{v}^{2}} \right),
  \end{aligned}\\
  G_0 = &
  \begin{aligned}[t]
  &1+256\text{ }{{u}^{2}}{{\left( 1-u \right)}^{2}}{{\left( 1-2v \right)}^{2}},
  \end{aligned}\\
  e_0 = &
  \begin{aligned}[t]
  &-32v\left( 1-v \right),
   \end{aligned}\\
  f_0 = &
  \begin{aligned}[t]
  &16\left( 1-2u \right)\left( 1-2v \right),
  \end{aligned}\\
  g_0 = &
  \begin{aligned}[t]
  &-32\text{ }u\text{ }(1-u).\label{humpg0}
  \end{aligned}
\end{align}
Thus, eq.~ \eqref{mcn} for $n=0$, along with eqs.~\eqref{humpE0} to \eqref{humpg0} gives
\begin{equation}\label{mc0s1}
  \begin{split}
    H_0 (u,v)  = &  32 (v (-1+v) \left(1+256 u^2 (-1+u)^2  (1-2 v)^2\right)  -256 \,u\, v  (-1+u) (-1+v) \\& \,  (1-2 u)^2  (1-2 v)^2   +(-1+u) u \left(1+256 v^2 (1-2 u)^2 (-1+v)^2 \right)),
  \end{split}
\end{equation}
as shown in  Fig.~\ref{mc0surface1}.
\noindent Substituting the value of $H_0$ from eq.~\eqref{mc0s1} along with eq.~\eqref{humpbuv} in  eqs.~ \eqref{VSN}~\eqref{vpn} for $n=0$ gives us the following expression
\begin{equation}\label{x1uvsurface1}
  \begin{split}
    \mathbf{x_1}(u,v,t) & = (u,v,16 (1-u) u (1-v) v + t (1-u) u (1-v) v (-2 (16-32 u-32 v+64 u v)  \\
      & \qquad (16 (1-u) u (1-v) -16 (1-u) u v) (16 (1-u) (1-v) v-16 u (1-v) v)  + \\
      & \qquad (-32 v+32 v^2) (1+(16 (1-u) u (1-v) -16 (1-u) u v)^2)+(-32 u+32 u^2) \\
      & \qquad (1+(16 (1-u) (1-v) v-16 u (1-v) v)^2))).
  \end{split}
\end{equation}
Denoting the fundamental magnitudes for this variational surface by  $E_{1}(u,v,t)$, $F_{1}(u,v,t)$, $G_{1}(u,v,t)$, $e_{1}(u,v,t)$, $f_{1}(u,v,t)$, $g_{1}(u,v,t)$  and plugging  in these  values in eq.~ \eqref{MCN} for $n=1$, we find the expression for $H_{1}(u,v,t)$, mean curvature of $\mathbf{x_1} (u,v,t)$. These fundamental magnitudes, the mean curvature  $H_{1}(u,v,t)$ and the coefficients $p_{i}(u,v)$ of $t^i  \text{\hspace{0.1cm} for \hspace{0.1cm}}   i = 0,1,2,3 $ in the expansion of $ H_{1}(u,v,t)$ are given in Appendix ~\ref{appendixa}. $ H_{1}(u,v,t)$ is a polynomial in $t$ and thus  $ H_{1}^{2}(u,v,t)$ is polynomial in $t$ as well. We find the non-zero coefficients $q_{i}(u,v)$ of $t^i \hspace{0.1cm}\text{for} \hspace{0.1cm} i = 0,1,2,3,4,5,6$ and integrate these coefficients for $0 \leq u,v \leq 1$ as mentioned in eq.~ \eqref{rmsn}  to get an
expression for the mean square of the mean curvature \eqref{rmsn} for $n=1$, as a polynomial in $t$. This is
\begin{equation}\label{h1nsq}
  \mu_{1}^{2} (t) = 1637.65 - 20425 t+195725 t^2-898809 t^3+ 2.98414 \times 10^6 t^4-5.10679 \times 10^6 t^5 + 4.1912 \times 10^6 t^6,
\end{equation}
shown in Fig.~\ref{graphh1nsqt}.
Minimizing this polynomial for $t$ gives  us  $t_{min1} = 0.088933$. We find the variationally improved surface $\mathbf x_{1}\left(u,v \right)$ eq.~ \eqref{VS1} for this minimum value of $t$, that is,
\begin{equation}\label{x1uvsurface1tmin}
\begin{split}
    \mathbf{x_{1}}(u,v) & = (u,v,(-1+u) u (-1+v) v (16+v (-2.84585+2.84585 v)+u^4 v (2185.61+v (-8013.91 \\&\qquad +(11656.6-5828.3 v) v))+u^3 v (-4371.22+v (16027.8+v  (-23313.2+11656.6 v))) \\&\qquad +u^2 (2.84585+v (2914.15+v (-10928.1+  (16027.8-8013.91 v) v)))+u (-2.84585 \\&\qquad +v (-728.537+v (2914.15+v (-4371.22+2185.61 v)))))),
\end{split}
\end{equation}
shown in Fig. ~\ref{graphvarx1Mint}. For this  $t_{min1}$ the mean curvature of $\mathbf{x_1}(u,v)$ is shown in Fig.~\ref{graphcH1Mint}. The initial area  of the  surface~$\mathbf{x_{0}}(u,v)$ (using eq.~\eqref{area0}) is  $2.494519$ and that of the surface~\eqref{x1uvsurface1tmin} (using eq.~\eqref{arean}) is  $2.11589$  for $t_{min1} = 0.0889$. The percentage decrease in the original area in this case is  $p_{01}=15.1784$. Substituting $H_1(u,v)$ (shown in Figure~\ref{graphcH1Mint}) in eq.~\eqref{vpn} along with eq.~\eqref{humpbuv}  for $n=1$  in eq.~ \eqref{VSN} results in an expression for the variational surface  $\mathbf{x_{2}}(u,v,t)$. We find  the fundamental magnitudes $E_{2}(u,v,t)$, $F_{2}(u,v,t)$, $G_{2}(u,v,t)$, $e_{2}(u,v,t)$, $f_{2}(u,v,t)$, $g_{2}(u,v,t)$ for this variational surface $\mathbf{x_{2}}(u,v,t)$ and insert these fundamental
 magnitudes in eq.~ \eqref{MCN} to get the expression for $H_{2}(u,v,t)$, mean curvature of $\mathbf{x_2} (u,v)$, and thus mean square of mean curvature $ H_{2}(u,v,t)$ (eq.~\eqref{rmsn} for $n=2$) gives the following expression,
\begin{equation}\label{h2nsq}
    \mu_{2}^{2} (t) =897.323 -14022.3 t + 207068 t^2-1.0771\times 10^6 t^3+6.4546 \times 10^6 t^4-9.9155\times 10^6 t^5+1.9927\times 10^7 t^6.
\end{equation}
Minimizing this expression results in  $t_{min2} = 0.0441$. For this $t_{min2} = 0.0441$ we find the variationally improved surface $\mathbf{x_{2}}(u,v)$  shown in Fig.~\eqref{graphvarx2Mint2}. The  numerator of the mean curvature of $\mathbf{x_{2}}(u,v)$ is shown in  Fig.~\ref{graphcH2Mint2} and the quantity $m_2(u,v)$ is shown in Fig.~\ref{m2graph}. A summary of related results is provided in table~\ref{table2}.
\subsection{The Bilinear Interpolant \eqref{eqbinterpolation1} Spanned by Four Boundary Non-Coplanar Straight Lines}\label{hessfalbis}
\begin{figure}[htb!]
\begin{center}
\includegraphics[width=60mm]{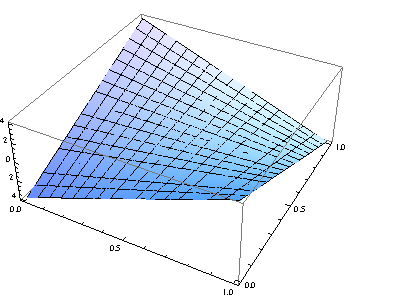}
\end{center}
\caption{$H_0(u,v)$, the numerator of the  mean curvature of the initial surface $\mathbf{x}_0(u,v)$, a bilinear interpolant for $0\leq u,v \leq 1$.} \label{mc0surface2}
\end{figure}

\begin{figure}[htb!]
\begin{center}
\includegraphics[width=60mm]{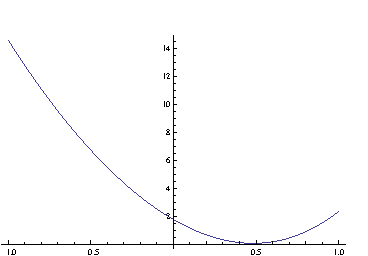}
\end{center}
\caption{Mean square $\mu_{1}^{2} (t)$ of the mean curvature of the surface $\mathbf{x}_1(u,v,t)$.} \label{bih1nsq}
\end{figure}

In the previous  subsection we have applied the algorithm eq.~\eqref{VSN} along with eq.~\eqref{vpn} (section \ref{basics}) to a surface for which the corresponding minimal surface is known. We have judged that our algorithm can significantly decrease area where the area can be decreased.  Below, we apply this algorithm to an important class of surfaces, namely the bilinear interpolant  where the corresponding minimal  area is not explicitly known. Specifically, we have taken the initial  surface $\mathbf{x}_{0}(u,v)$ given by eq.~ \eqref{eqbinterpolation1} for $0\leq u,v\leq 1$. It is  bounded by four non-coplanar straight lines. A convenient  alternative for the  unit normal $\mathbf N(u,v)$ in  eq.~\eqref{VSN} is $\mathbf k=(-1,0,0)$, making a small angle with the  unit normal $\mathbf N(u,v)$  to   the surface given by eq.~ \eqref{eqbinterpolation1}. The corresponding function $b(u,v)$ assuring  that the  variation  of the surface at the boundary curves at $u=0, u=1, v=0, v=1$ vanishes is given by
\begin{equation}
  b(u,v)=u v(1-u)(1-v).
\end{equation}
The surface given by \eqref{eqbinterpolation1} for $0 \leq u,v\leq 1$ is a non-minimal surface spanned by a boundary composed of non-coplanar straight lines.  The  fundamental magnitudes of this initial surface are
\begin{equation}\label{bifm}
    E_0=1+(1-2 v)^2, F_0=(1-2 u) (1-2 v), G_0=1+(1-2 u)^2, e_0=0, f_0=2, g_0=0.
\end{equation}
For $n=0$, eq.~ \eqref{mcn} along with eq.~  \eqref{bifm} gives
\begin{equation}\label{mc0s2}
    H_{0}=-4 (1-2 u) (1-2 v),
\end{equation}
as shown in the Fig.~\ref{mc0surface2}.
The variationally improved surface eq.~ \eqref{VSN} along with its fundamental magnitudes as functions of $u,v$ and $t$ are
\begin{align}
    \mathbf{x}_1(u,v,t)=&
    \begin{aligned}[t]
    & ( u+v-2 u v + 4u v(1-u)(1-v)(1-2 u)(1-2v) t,v,u ),
    \end{aligned}\\
    E_1(u,v,t)=&
    \begin{aligned}[t]
    &1+(1-2 v)^2 \left(1-4 t \left(1-6 u+6 u^2\right) (-1+v) v\right)^2,\label{bifmE1}
    \end{aligned}\\
    F_1(u,v,t)=&
    \begin{aligned}[t]
    &(1-2 u) (1-2 v) \left(4 t \,(v^2-v)\left(6 u^2-6 u+1\right)  -1\right) \left(4 t  (u^2-u) \left(6 v^2-6 v+1\right)-1\right),
    \end{aligned}\\
    G_1(u,v,t)=&
    \begin{aligned}[t]
    &(1-2 u)^2 \left(1-4 t (u-1) u \left(6 v^2-6 v+1\right)\right)^2+1,
    \end{aligned}\\
    e_1(u,vt)=&
    \begin{aligned}[t]
    &-24 t (-1+2 u) v \left(1-3 v+2 v^2\right),
    \end{aligned}\\
    f_1(u,v,t)=&
    \begin{aligned}[t]
    &-2 \left(-1+2 t \left(1-6 u+6 u^2\right) \left(1-6 v+6 v^2\right)\right),
    \end{aligned}\\
    g_1(u,v,t)=&
    \begin{aligned}[t]
    &-24 t u \left(1-3 u+2 u^2\right) (-1+2 v).\label{bifmg1}
    \end{aligned}
\end{align}
Thus, eq.~\eqref{MCN} along with eqs.~\eqref{bifmE1} to \eqref{bifmg1} gives an expression for $H_1(u,v,t)$ and hence the coefficients $p_i(u,v)$ of $t^i$ for $i=0, 1, 2, 3$ in the expansion of $H_1(u,v,t)$ as mentioned in the eq.~\eqref{MCN}, which are given in Appendix ~\ref{appendixb}. $H_{1}(u,v,t)$ is a polynomial in $t$ and thus  $ H_{1}^{2}(u,v,t)$ is polynomial in $t$ as well. We find the non-zero coefficients $q_{i}(u,v)$ of $t^i \hspace{0.1cm}\text{for} \hspace{0.1cm} i = 0,1,2,3,4,5,6$ and integrate these coefficients for $0 \leq u,v \leq 1$, as mentioned in eq.~ \eqref{rmsn},  to get an expression for the mean square of mean curvature as a polynomial in $t$, given by
\begin{equation}\label{mc0bi}
      \mu_{1}^{2} (t) = 1.7778 -6.8267 t+6.4261 t^2+0.6966  t^3+0.2648 t^4+0.0076 t^5+0.0009 t^6.
\end{equation}
Minimizing this  polynomial for $t$ gives us $t_{min}=0.4836$. The mean square of mean curvature of the surface $\mathbf{x}_1(u,v)$ as a function of $t$ is shown in Fig.~\ref{bih1nsq}.

\begin{figure}[htb!]
\begin{center}
\includegraphics[width=60mm]{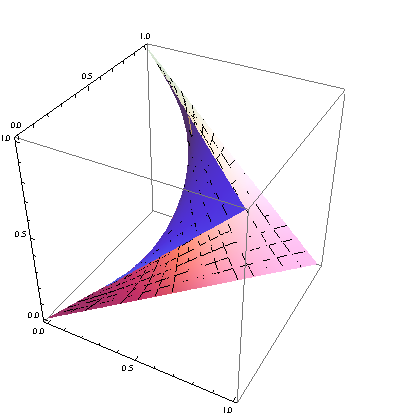}
\end{center}
\caption{$\mathbf{x}_1(u,v,t)$ for $t=t_{min}=0.4836 $ for $0\leq u,v\leq 1$} \label{bivarx1Mint}
\end{figure}

\begin{figure}[htb!]
\begin{center}
\includegraphics[width=60mm]{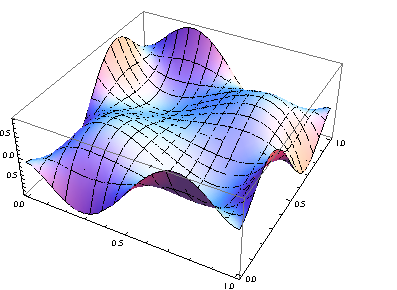}
\end{center}
\caption{ Mean curvature of the surface $\mathbf{x_1}(u,v)$ for $t= t_{min} = 0.4836$ for $0\leq u,v\leq 1$} \label{bicH1Mint}
\end{figure}

\begin{figure}[htb!]
\begin{center}
\includegraphics[width=60mm]{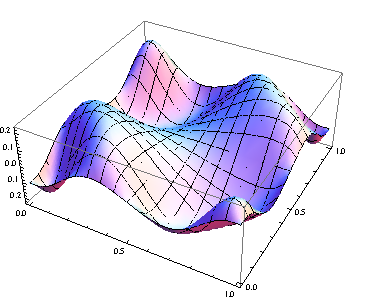}
\end{center}
\caption{Mean curvature of the surface $\mathbf{x}_2(u,v)$ for $t_{min2}=0.0881$.} \label{bicH2Mint}
\end{figure}

\noindent We find the variationally improved surface eq.~\eqref{VS1} for this $t_{min}=0.4836$,
\begin{equation}\label{x1uvbi}
    \mathbf{x}_1(u,v)=\left(u+v-2 u v+1.9345\, u\, v  (1- u) (1- v)  (1- 2 u) (1- 2 v) , v, u\right),
\end{equation}
shown in Fig.~\ref{bivarx1Mint}.
The mean curvature of $\mathbf{x_1}(u,v)$ is shown in Fig.~\ref{bicH1Mint}.
The initial area of the surface $\mathbf{x}_0(u,v)$ (eq.~\eqref{eqbinterpolation1}) (using eq.~\eqref{area0}) is  $A_0 = 1.2808$ and that of surface \eqref{x1uvbi} is $A_1 = 1.2793$. In a similar way, we find the second order variation of the bilinear interpolant whose mean square $ms$ of mean curvature, according to eq.~\eqref{rmsn}, is
\begin{equation}\label{bismc}
\mu_2^2(t)= 0.0728 -1.2952 t+7.3524 t^2-0.0058  t^3+0.0147 t^4-0.0001 t^5+0.00003 t^6,
\end{equation}
as shown in  Fig.~\ref{bih2nsqt} as a function of $t$. (Details are given in Appendix ~\ref{appendixb}). Here  $t_{min2}=0.0881$ for which the variational surface  $\mathbf{x}_2(u,v)$ is   shown in  Fig.~\ref{bivarx2Mint} and  given in Appendix ~\ref{appendixb}. The  curvature of this surface is  shown in Fig.~\ref{bicH2Mint}.
\begin{figure}[htb!]
\begin{center}
\includegraphics[width=60mm]{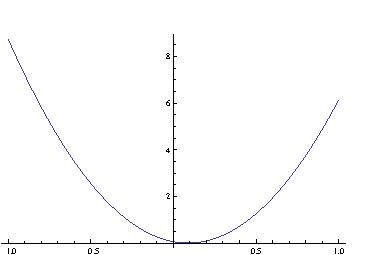}
\end{center}
\caption{Mean square $\mu_2^2(t)$ of the mean curvature  of the bilinear interpolant  $\mathbf{x}_2(u,v,t)$ as a function of $t$.} \label{bih2nsqt}
\end{figure}

\begin{figure}[htb!]
\begin{center}
\includegraphics[width=60mm]{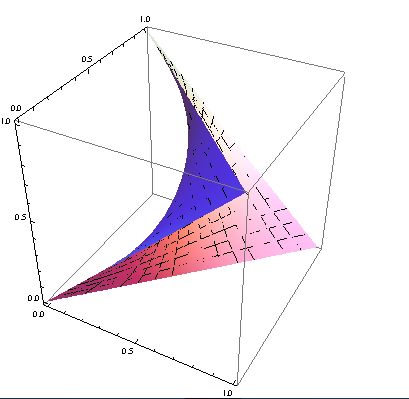}
\end{center}
\caption{$\mathbf{x}_2(u,v,t)$ for $t=t_{min2}=0.0881$ for $0\leq u,v\leq 1$.} \label{bivarx2Mint}
\end{figure}
\begin{figure}[htb!]
\begin{center}
\includegraphics[width=60mm]{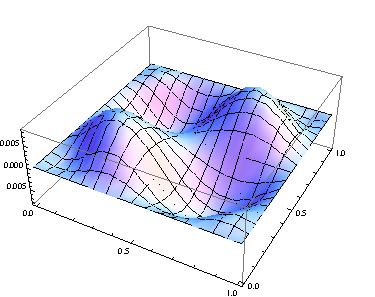}
\end{center}
\caption{$m_2(u,v)$ (eq.~\eqref{vpn} for $n=2$) of the surface $\mathbf{x_3}(u,v)$ for $0\leq u,v\leq 1$.} \label{bim2graph}
\end{figure}

\noindent The initial area (calculated using eq.~\eqref{arean}) of surface ~$\mathbf{x_{0}}(u,v)$  of eq.~ \eqref{eqbinterpolation1} is $1.280789$  and that (eq.~\eqref{arean} for $n=1$) of surface ~$\mathbf{x_{1}}(u,v)$ of eq.~\eqref{x1uvbi}  is $1.27936$  for $t_{min} =0.4836$. The area (eq.~\eqref{arean} for $n=2$) of surface ~$\mathbf{x_{2}}(u,v)$ of eq.~\eqref {vcpr2bi}  is $1.279301$. The percentage decrease in area with respect to that of eq.~\eqref{x1uvbi} is  0.004382. We have not been able to find higher order variational surfaces for $n\geq 3$ as in this case our computer program becomes unresponsive for higher iterations. However to foresee that  a further reduction is possible, the variational quantity $m_2(u,v)$ for ~$\mathbf{x_{3}}(u,v)$ is shown in Fig.~\ref{bim2graph}. The algorithm reduces area less significantly  for the bilinear interpolant.  We suggest this relative stability indicates that the bilinear interpolant is  a quasi-minimal surface. The ratio of $rms$ of Gaussian curvature to the $ms$ of the mean curvature obtained for successive surfaces decreases for  the hemiellipsoid,  whereas this ratio increases for the hump-like surface and the bilinear interpolant.   Moreover, it can be seen from table \ref{table2} that  by using the actual $H_0$, the percentage decrease in area of the hemiellipsoid is slightly more than the decrease obtained by replacing $H_0$ by 1. A  summary of these results is presented in table~\ref{table2}.
\begin{center}
\begin{table}[htb!]
\caption{\textbf{Summary of Results for the Surfaces}}
\smallskip
\begin{quote}
  A summary of results for the surfaces $\mathbf{x}_{i} (u,v)$, with starting surface written in the first column. $A_i$ are the decreasing areas, $p_{ij}$ and $q_{ij}$ ($i<j$) are the percentage decreases in areas. $\nu_{i}$  are the $rms$ of Gaussian curvature and $\mu_{i} $ are the $rms$ of the mean curvature.  The last column reports the optimal value $t_{min}$ of our variational parameter $t$ for each of the cases.
\end{quote}
\vspace{3mm}
\begin{tabular}{|p{3.8cm} |p{1.5cm}|p{1.5cm}|p{1.5cm}|p{1.5cm}|p{1.5cm}| }
\hline
\multicolumn{6}{ |c| }{\textbf{Hemiellipsoid, Hump-Like Surface and Bilinear Interpolant}} \\
\hline
$\mathbf{x}_{i} (u,v)$  &  $ A_{i}$ &  $p_{ij}$ & $q_{ij}$ & $\nu_{i}/ \mu_{i}^{2}$ & $t_{min}$  \\
\hline

\multirow{3}{*}{Hemiellipsoid   ($H_0=1$)} &  6.28319  &  -  &  -  & 0.38985  &   \\
 &  4.70625   &  50.1954  &  -   & 0.01369 &  -0.35157\\
 &  4.4403  & 8.46711      & -    & 0.01308 &  0.07064  \\
 & 4.4025      &  1.20025    & -   & 0.01307 &  0.02264  \\
  \hline
 \multirow{3}{*}{Hemiellipsoid} &  6.28319  &  -  &  -  & 0.38985  &  - \\
 &  4.32641   &  62.2861 &  -   &  0.19644 & 0.14825  \\
 & 4.23731  &  2.83624             &  -       &  0.18095 &  0.02297 \\
 \hline
 \multirow{3}{*}{Hump-Like Surface} & 2.49452 & - & -  & 0.03456 & - \\
 & 2.11589 & 25.3341 & - &  0.05207  &   0.08893 \\
 & 1.92788 &12.58 & -  &  0.09449  & 0.04409\\
 \hline
\multirow{2}{*}{Bilinear Interpolant} &  1.280789 &  -  &  -  &  2.25  &  - \\
&   1.27936  &  - & 0.11157 & 58.9258 & 0.48364 \\
& 1.2793 &  - & 0.00461  & 272.152 &  0.0881 \\
\hline
\end{tabular} \label{table2}
\end{table}
\end{center}
\pagebreak
\section{Conclusions} \label{conclusion}
We have discussed how to reduce the surface area of the non-minimal surfaces given by eqs.~\eqref{hemiellipsoid}, \eqref{eqsurface1} and \eqref{eqbinterpolation1} by using a variational technique with an appropriate number of iterations to improve the surfaces $\mathbf x_{i}(u,v)$.  The first two cases are meant to test our curvature algorithm and the third one applies this algorithm to a surface for which there is no corresponding known minimal surface.   The percentage decrease in area of the hemiellipsoid eq.~\eqref{hemiellipsoid} is given for two cases: one for which $H_0$ is replaced by 1 and one for which $H_0$ is the usual numerator  of the mean curvature defined by eq.~\eqref{mcn}. In the notation of eq.~\eqref{percentdecknown}, the area reductions in the former case are $p_ {01}= 50.1954, p_ {02}=58.6625$ and  $p_ {03}=59.8628$, and in the latter case are $p_ {01}= 62.2861$ and $p_ {02}=65.1223$. For the hump-like surface eq. \eqref{eqsurface1}, the percentage decreases in area are $p_ {01}= 25.3341$ and $p_ {02}=37.9141$.  This indicates that the algorithm attains at least a local minimum of the area in these two cases. The percentage  decrease (resulting through the corresponding slightly different definition given by eq.~\eqref{percentdecunknown}) in area of the surface eq.~\eqref{eqbinterpolation1} is much less i.e. $q_ {01}= 0.1116$  and $q_ {02}=0.1162$. This means that in the case of bilinear interpolant with  unknown minimal surface, the variational improvement does \emph{not} result in a significant decrease in the area. This  saturation indicates that the bilinear interpolant is at least a local minima.   Thus the bilinear interpolant is indicated as (or may well be) a critical point of the area.
\appendix

\section{Expressions Used in Section~\ref{hessfalhs}}\label{appendixa}
Here   $E_{1}(u,v,t)$, $F_{1}(u,v,t)$, $G_{1}(u,v,t)$, $e_{1}(u,v,t)$, $f_{1}(u,v,t)$, $g_{1}(u,v,t)$  denote  the  fundamental magnitudes of the first variational surface $\mathbf{x}_1(u,v,t)$ obtained by plugging  eq.~\eqref{eqsurface1} in eq.~\eqref{VSN} for $n=0$  where $p_i$  denotes the coefficient  of  $t^i$ ($i=0,1,2,3$) in the polynomial expansion of $H_n$ (eq.~\eqref{MCN}) for $n=1$. $H_1(u,v,t)$ is the numerator of the mean curvature of the surface $\mathbf{x}_1(u,v,t)$.
\begin{align}
    E_{1}(u,v,t)  =&
    \begin{aligned}[t]
    & 1 + 256 (1 - 2 u)^2 (-1 + v)^2 v^2 (1 +    2 t ((-1 + v) v + 1536 u^3 (-1 + v) v (3 + 8 (-1 + v) v) \\& -
       768 u^4 (-1 + v) v (3 + 8 (-1 + v) v) + 2 u (-1 + 256 (-1 + v) v (1 + 3 (-1 + v) v)) \\& +
       2 u^2 (1 - 128 (-1 + v) v (11 + 30 (-1 + v) v))))^2,\label{E1(u,v,t)}
    \end{aligned}\\
     F_{1}(u,v,t) =&
     \begin{aligned}[t]
     & 256 (-1+u) u (-1+2 u) (-1+v) v (-1+2 v) (-1+2 t (-2 (-1+u) u+v-256 (-1+u) \\& u  (-2+3 u)  (-1+3 u) v+(-1+256 (-1+u) u (8+33 (-1+u) u)) v^2-3072 (1-2 u)^2 \\& (-1+u) u v^3  + 1536 (1-2 u)^2 (-1+u) u v^4)) (-1+2  t (-2 (-1+v) v-3072 u^3 \\& (1-2 v)^2 (-1+v) v  + 1536 u^4 (1-2 v)^2 (-1+v) v+u (1-256 (-1+v) v \\& (-2+3 v) (-1+3 v))+u^2 (-1+256 (-1+v) v (8+33 (-1+v) v)))),\label{F1(u,v,t)}
     \end{aligned}\\
     G_{1}(u,v,t) =&
    \begin{aligned}[t]
    & 1+256 (-1+u)^2 u^2 (1-2 v)^2 (1+2 t ((-1+u) u-2 v+512 (-1+u)  u (1+3 (-1+u) u) v \\& +2 (1-128 (-1+u) u (11+30 (-1+u) u)) v^2 +1536 (-1+u) u (3+8 (-1+u) u) v^3-768 \\&  (-1+u) u (3+8 (-1+u) u) v^4))^2,\label{G1(u,v,t)}
    \end{aligned}\\
      e_{1}(u,v,t) =&
     \begin{aligned}[t]
     &  -32 (-1+v) v (-1+2 t (-1+6 u-6 u^2-15 (17+256 (-1+u) u (1+3 (-1+u) u)) v \\& +3 (341+256 (-1+u) u (19+55 (-1+u) u)) v^2-1536 (1-2 u)^2 (1+10 (-1+u) u) v^3 \\& +768 (1-2 u)^2 (1+10 (-1+u) u) v^4)),\label{e1(u,v,t)}
     \end{aligned}\\
     f_{1}(u,v,t)  =&
    \begin{aligned}[t]
    & -16 (-1+2 u) (-1+2 v) (-1+4 t (-(-1+u) u+v-256 (-1+u) u (-2+3 u) (-1+3 u) v \\& + (-1+256 (-1+u) u (11+45 (-1+u) u)) v^2-4608 (1-2 u)^2 (-1+u) u v^3 \\& +2304 (1-2 u)^2 (-1+u) u v^4)),\label{f1(u,v,t)}
    \end{aligned}\\
     g_1(u,v,t) =&
    \begin{aligned}[t]
    & -32 (-1+u) u (-1+2 t (-1-6 (-1+v) v-1536 u^3 (1-2 v)^2 (1+10 (-1+v) v)  \\
      & + 768 u^4 (1-2 v)^2 (1+10 (-1+v) v)-15 u (17+256 (-1+v) v (1+3 (-1+v) v))  \\& +3 u^2 (341+256 (-1+v) v (19+55 (-1+v) v)))),\label{g1(u,v,t)}
    \end{aligned}
  \end{align}

\begin{equation}
  \begin{split}
 & H_{1} (u,v,t)  =\, \\&  -32 (-1+v) v (-1+2 t (-1+6 u-6 u^2-15 (17+256 (-1+u) u (1+3 (-1+u) u)) v+3 (341 \, \\& + 256 (-1+u)
u (19+55 (-1+u) u)) v^2-1536 (1-2 u)^2 (1+10 (-1+u) u) v^3+768 (1-2 u)^2 (1+10 \, \\& (-1+u) u)
v^4)) (1+256 (-1+u)^2 u^2 (1-2 v)^2 (1+2 t ((-1+u) u-2 v+512 (-1+u) u (1+3 (-1+u) u) v\, \\& +2
(1-128 (-1+u) u (11+30 (-1+u) u)) v^2+1536 (-1+u) u (3+8 (-1+u) u) v^3-768 (-1+u) u \, \\&
(3+8 (-1+u) u) v^4))^2)+8192 (1-2 u)^2 (-1+u) u (1-2 v)^2 (-1+v) v (-1+2 t (-2 (-1+u) u+v \, \\&
-256 (-1+u) u (-2+3 u) (-1+3 u) v+(-1+256 (-1+u) u (8+33 (-1+u) u)) v^2-3072 (1-2 u)^2  \, \\&
(-1+u) u v^3+1536 (1-2 u)^2 (-1+u) u v^4)) (-1+4 t (-(-1+u) u+v-256 (-1+u) u (-2+3 u) (-1+3 u) \, \\& v
+(-1+256 (-1+u) u (11+45 (-1+u) u)) v^2-4608 (1-2 u)^2 (-1+u) u v^3+2304 (1-2 u)^2 (-1+u) \, \\&
u v^4)) (-1+2 t (-2 (-1+v) v-3072 u^3 (1-2 v)^2 (-1+v) v+1536 u^4 (1-2 v)^2 (-1+v) v+u (1-256 \, \\&
(-1+v) v (-2+3 v) (-1+3 v))+u^2 (-1+256 (-1+v) v (8+33 (-1+v) v))))-32 (-1+u) u (-1+2 t \, \\& (-1-6
(-1+v) v-1536 u^3 (1-2 v)^2 (1+10 (-1+v) v)+768 u^4 (1-2 v)^2 (1+10 (-1+v) v)-15 u (17 \, \\&
+256 (-1+v) v (1+3 (-1+v) v))+3 u^2 (341+256 (-1+v) v (19+55 (-1+v) v)))) (1+256 (1-2 u)^2 \, \\&
(-1+v)^2 v^2 (1+2 t ((-1+v) v+1536 u^3 (-1+v) v (3+8 (-1+v) v)-768 u^4 (-1+v) v (3+8 (-1+v) v) \, \\&
+2 u (-1+256 (-1+v) v (1+3 (-1+v) v))+2 u^2 (1-128 (-1+v) v (11+30 (-1+v) v))))^2).\label{H1(u,v,t)}
 \end{split}
 \end{equation}

\begin{equation}\label{h0}
\begin{split}
  p_{0} (u,v) & = -32 (-(-1+u) u+v-256 (-1+u) u (1+3 (-1+u) u) v+(-1+256 (-1+u) u (4+11 \\&\qquad  (-1+u) u)) v^2  -512 (-1+u) u (3+8 (-1+u) u) v^3+256 (-1+u) u (3+8 (-1+u) u) v^4),
\end{split}
\end{equation}

\begin{equation}\label{h1}
\begin{split}
  p_{1} (u,v) & = 64 (-786432 u^7 (-1+v)^2 v^2 (21+4 (-1+v) v (31+48 (-1+v) v))+196608 u^8 (-1+v)^2 \\&\qquad v^2 (21+4 (-1+v) v (31+48 (-1+v) v))-(-1+v) v (-1+3 (-1+v) v (85+256 (-1+v) v))\\&\qquad + u (-1+4 (-1+v) v (-3+64 (-1+v) v (19+54 (-1+v) v)))-768 u^5 (-3+2 (-1+v) v \\&\qquad (-27+2 (-1+v) v (7835+256 (-1+v) v (184+285 (-1+v) v))))+256 u^6 (-3+2 (-1+v) v \\&\qquad (-27+2 (-1+v) v (26651+768 (-1+v) v (206+319 (-1+v) v))))+2 u^2 (-127+2 (-1+v) v \\&\qquad (-1213+64 (-1+v) v (1405+2 (-1+v) v (4891+8064 (-1+v) v))))-2 u^3 (-639+512 (-1+v)\\&\qquad  v (-23+(-1+v) v (3171+2 (-1+v) v (10139+15936 (-1+v) v))))+u^4 (-2559+256 (-1+v)\\&\qquad  v (-181+4 (-1+v) v (11765+(-1+v) v (72347+112320 (-1+v) v))))),
\end{split}
\end{equation}

\begin{equation}\label{h2}
\begin{split}
  p_2 (u,v) & = -32768 (-1+u) u (-1+v) v ((-1+v)^2 v^2 (6+19 (-1+v) v)-11796480 u^9 (-1+v)^2 v^2 (14+(-1+v) \\&\qquad  v (119+32 (-1+v) v (11+12 (-1+v) v)))+2359296 u^{10} (-1+v)^2 v^2 (14+(-1+v) v (119+32  \\&\qquad (-1+v) v (11+12 (-1+v) v)))+u (-1+v) v (-9+2 (-1+v) v (2029+2 (-1+v) v (7093+12672  \\&\qquad (-1+v) v)))-6144 u^7 (-1+v) v (-33+16 (-1+v) v (4453+6 (-1+v) v (6383+128 (-1+v) v  \\&\qquad (149+164 (-1+v) v))))+1536 u^8 (-1+v) v (-33+16 (-1+v) v (14533+6 (-1+v) v (20663+128  \\&\qquad (-1+v) v (479+524 (-1+v) v))))-u^3 (31+8 (-1+v) v (-4561+8 (-1+v) v (143089+32  \\&\qquad (-1+v) v (43313+48 (-1+v) v (2981+3528 (-1+v) v)))))+u^2 (6+(-1+v) v (-4049+2  \\&\qquad (-1+v) v (367651+2 (-1+v) v (1912651+384 (-1+v) v (17455+21504 (-1+v) v)))))-3 u^5  \\&\qquad (19+4 (-1+v) v (-23989+256 (-1+v) v (53227+12 (-1+v) v (39553+4 (-1+v) v (30469  \\&\qquad + 34384 (-1+v) v)))))+u^6 (19+4 (-1+v) v (-83125+256 (-1+v) v (326171+12 (-1+v) v  \\&\qquad (237065+12 (-1+v) v (59735+66416 (-1+v) v)))))+u^4 (63+2 (-1+v) v (-69931+16  \\&\qquad (-1+v) v (1582433+32 (-1+v) v (456095+24 (-1+v) v (60259+69408 (-1+v) v)))))),
\end{split}
\end{equation}
\begin{equation}\label{h3}
\begin{split}
p_{3}(u,v) &  = 65536 (-1+u) u (-1+v) v (-(-1+v)^3 v^3 (3+10 (-1+v) v)-25367150592 u^{13} (1-2 v)^2 \\&\qquad  (-1+v)^3 v^3 (3+8 (-1+v) v) (2+(-1+v) v (9+16 (-1+v) v))+3623878656 u^{14} (1-2 v)^2 \\&\qquad  (-1+v)^3 v^3 (3+8 (-1+v) v) (2+(-1+v) v (9+16 (-1+v) v))-u (-1+v)^2 v^2 (-6+(-1+v) v (2793 \\&\qquad  +8 (-1+v) v (2525+4704 (-1+v) v)))+2359296 u^{12} (-1+v)^2 v^2 (-21+(-1+v) v (203161 \\&\qquad  +64 (-1+v) v (35519+(-1+v) v (155053+768 (-1+v) v (421+355 (-1+v) v)))))-14155776 u^{11} \\&\qquad   (-1+v)^2 v^2 (-21+(-1+v) v (63385+64 (-1+v) v (11131+(-1+v) v (48765+128 (-1+v) v  \\&\qquad   (797+674 (-1+v) v)))))-u^2 (-1+v) v (6+(-1+v) v (-5090+(-1+v) v (999971+8 (-1+v) v  \\&\qquad   (1355491+96 (-1+v) v (51427+64512 (-1+v) v)))))-768 u^9 (-1+v) v (245+4 (-1+v) v \\&\qquad   (-386905+8 (-1+v) v (38922251+256 (-1+v) v (1740737+24 (-1+v) v (322558+(-1+v) v (683777 \\&\qquad   +584928 (-1+v) v))))))+768 u^{10} (-1+v) v (49+4 (-1+v) v (-254789+8 (-1+v) v (45198319+256 \\&\qquad   (-1+v) v (1999045+24 (-1+v) v (367093+(-1+v) v (772045+655968 (-1+v) v))))))+u^3 (3-(-1+v) v \\&\qquad   (2781+4 (-1+v) v (-248143+2 (-1+v) v (14160189+128 (-1+v) v (1627565+6 (-1+v) v (1477367 \\&\qquad   +384 (-1+v) v (9583+9216 (-1+v) v)))))))-8 u^7 (-5+4 (-1+v) v (14285+128 (-1+v) v (-179703+4  \\&\qquad   (-1+v) v (15643685+32 (-1+v) v (5816093+48 (-1+v) v (556153+36 (-1+v) v (33637+29392  \\&\qquad   (-1+v) v)))))))+u^4 (-19+(-1+v) v (28573+4 (-1+v) v (-3466831+(-1+v) v (508303293+128 (-1+v) v  \\&\qquad   (53765501+18 (-1+v) v (15240289+512 (-1+v) v (70459+65280 (-1+v) v)))))))+2 u^8 (-5+4 (-1+v) v  \\&\qquad   (49565+64 (-1+v) v (-2244879+64 (-1+v) v (17885581+4 (-1+v) v (52018853+48 (-1+v) v  \\&\qquad   (4883029+72 (-1+v) v (145395+125528 (-1+v) v)))))))+u^5 (49-(-1+v) v (126811+4 (-1+v) \\&\qquad   v (-21499593+256 (-1+v) v (15874307+8 (-1+v) v (25388918+3 (-1+v) v (41251309+288 \\&\qquad   (-1+v) v (327257+295936 (-1+v) v)))))))+u^6 (-63+(-1+v) v (312153+4 (-1+v) v \\&\qquad   (-78350147+256 (-1+v) v (76888369+8 (-1+v) v (117926386+3 (-1+v) v (185099599 \\&\qquad   +288 (-1+v) v (1428851+1267712 (-1+v) v)))))))).
\end{split}
\end{equation}
\newpage
\section{Expressions Used in Section~\ref{hessfalbis}} \label{appendixb}
Here $H_1(u,v,t)$ is the numerator of the mean curvature of the surface $\mathbf{x}_1(u,v,t)$ obtained by plugging eq.~\eqref{eqbinterpolation1} in eq.~\eqref{VSN} for $n=0$   where $p_i$  denotes the coefficient  of $t^i$ ($i=0,1,2,3$)  in the polynomial expansion of $H_n$ (eq.~\eqref{MCN}) for $n=1$. The variational surface  $\mathbf{x_{2}}(u,v,t)$ is obtained by plugging  eq.~\eqref{x1uvbi} in  eq.~ \eqref{VSN} for $n=1$. Minimizing the polynomial   in eq.~\eqref{bismc} for $t$ gives us  $t_{min}= 0.0881$, which results in the surface of lesser area given by $\mathbf{x_{2}}(u,v)$.
\begin{align}
    H_1(u,v,t) =&
    \begin{aligned}[t]
    & 4 (-6 t u (1-3 u+2 u^2) (-1+2 v) (1+(1-2 v)^2 (1-4 t (1-6 u+6 u^2) (-1+v) v)^2)+ (-1+2 u) \,\\& (-1+2 v)  (-1+4 t (1-6 u+6 u^2) (-1+v) v) (-1+4 t (-1+u) u (1-6 v+6 v^2)) (-1+2 t (1-6 u \,\\& +6 u^2) (1-6 v+6 v^2))-6 t (-1+2 u) v (1-3 v+2 v^2) (1+(1-2 u)^2 (1-4 t (-1+u)\,\\& u (1-6 v+6 v^2))^2)).
    \end{aligned}\\
    p_0(u,v)=&
    \begin{aligned}[t]
    &16 (1-2 u)^2 (1-2 v)^2,
    \end{aligned}\\
    p_1(u,v)=&
    \begin{aligned}[t]
    &-64 (1-2 u)^2 (1-2 v)^2 \left(1-2 v+2 v^2+u \left(-2+36 v-36 v^2\right)+u^2 \left(2-36 v+36 v^2\right)\right),
    \end{aligned}\\
    p_2(u,v) =&
    \begin{aligned}[t]
    & 64 (1-2 u)^2 (1-2 v)^2 (1-8 v+36 v^2-56 v^3+28 v^4-8 u (1-17 v+65 v^2-96 v^3+48 v^4) \,\\& -8 u^3 (7-96 v+636 v^2-1080 v^3+540 v^4)+4 u^4 (7-96 v+636 v^2-1080 v^3+540 v^4)+4 u^2  \,\\& (9-130 v+766 v^2-1272 v^3+636 v^4)),
    \end{aligned}\\
    p_3(u,v)=&
    \begin{aligned}[t]
    & -512 (1-2 u)^2 (1-2 v)^2 (v (-1+9 v-28 v^2+44 v^3-36 v^4+12 v^5)+u (-1+20 v-174 v^2 \,\\& +692 v^3-1306 v^4+1152 v^5-384 v^6)-36 u^5 (1-32 v+278 v^2-1212 v^3+2406 v^4  -2160 v^5 \,\\& +720 v^6)+12 u^6 (1-32 v+278 v^2-1212 v^3+2406 v^4-2160 v^5+720 v^6)+u^2 (9-174 v \,\\& +1432 v^2-5852 v^3+11266 v^4-10008 v^5+3336 v^6)-4 u^3 (7-173 v+1463 v^2\,\\&-6216 v^3+12198 v^4-10908 v^5+3636 v^6)+2 u^4 (22-653 v+5633 v^2-24396 v^3+\,\\&48288 v^4-43308 v^5+14436 v^6)).
    \end{aligned}
 \end{align}
\begin{equation}\label{vcpr2bi}
\begin{split}
    \mathbf{x}_2(u,v,t) & = (u+v-2 u v+1.9345 (-1+u) u (-1+2 u) (-1+v) v (-1+2 v)-t (1-u) u (1-v) v ( -2 ( 1  +\\
      & \quad ( -0.0655-11.6073 u+11.6073 u^2 ) v+( -5.8036+34.8218 u-34.8218 u^2 ) v^2+ ( 3.8691  -23.2145 u+\\
      & \quad 23.2145 u^2 ) v^3 ) ( 1  +u^2 ( -5.8036+34.8218 v-34.8218 v^2 )+u ( -0.06555-11.6073 v+11.6073 v^2 )+u^3\\
      & \quad  ( 3.8691  -23.2145 v+23.2145 v^2 ) ) ( 0.06546  +11.6073 v-11.6073 v^2+u^2 ( -11.6073+69.6435 v-\\
      & \quad 69.6435 v^2 )+u ( 11.6073  -69.6435 v+69.6435 v^2 ) )+u ( 11.6073  +u^2 (23.2145  -46.429 v)-23.2145 v+\\
      & \quad u (-34.8218+69.6435 v) ) ( 1+( 1  +( -0.0656-11.6073 u+11.6073 u^2 ) v+( -5.80363+34.8218 u-\\
      & \quad 34.8218 u^2 ) v^2+( 3.86908  -23.2145 u+23.2145 u^2 ) v^3 )^2 )+v ( 11.6073  -34.8218 v+23.2145 v^2+u \\
      & \quad ( -23.2145+69.6435 v-46.429 v^2 ) ) ( 1+( 1  +u^2 ( -5.8036+34.8218 v-34.8218 v^2 )+u ( -0.0655-\\
      & \quad 11.6073 v+11.6073 v^2 )+u^3 ( 3.8691  -23.2145 v+23.2145 v^2 ) )^2 ) ),v,u),
\end{split}
\end{equation}
\begin{equation}\label{bix2uv}
\begin{split}
  \mathbf{x}_2(u,v) & =(u (1-v)+(1-u) v+1.9345 (1-2 u) (1-u) u (1-2 v) (1-v) v-0.08807 (1-u) u (1-v) v (-2 (0.0655 +\\
      & \quad 11.6073 u-11.6073 u^2+11.6073 v-69.6435 u v+69.6435 u^2 v-11.6073 v^2+69.6435 u v^2-69.6435 u^2 v^2)\\
      & \quad  (1-2 u-0.4836 (-4 (1-2 u) (1-u) u (1-2 v) (1-v)+4 (1-2 u) (1-u) u (1-2 v) v+8 (1-2 u) (1-u) \\
      & \quad u (1-v) v)) (1-2 v-0.4836 (-4 (1-2 u) (1-u) (1-2 v) (1-v) v+4 (1-2 u) u (1-2 v) (1-v) v+8 (1-u)\\
      & \quad  u (1-2 v) (1-v) v))+(11.6073 v-23.2145 u v-34.8218 v^2+69.6435 u v^2+23.2145 v^3-46.429 u v^3) (1+\\
      & \quad (1-2 u-0.4836 (-4 (1-2 u) (1-u) u (1-2 v) (1-v)+4 (1-2 u) (1-u) u (1-2 v) v+8 (1-2 u) (1-u)\\
      & \quad  u (1-v) v))^2)+(11.6073 u-34.8218 u^2+23.2145 u^3-23.2145 u v+69.6435 u^2 v-46.429 u^3 v) (1+\\
      & \quad (1-2 v-0.4836 (-4 (1-2 u) (1-u) (1-2 v) (1-v) v+4 (1-2 u) u (1-2 v) (1-v) v+8 (1-u) u (1-2 v)\\
      & \quad  (1-v) v))^2)),(1-u) v+u v,u (1-v)+u v),
\end{split}
\end{equation}
\bibliographystyle{unsrt}
\bibliography{xbib}

\end{document}